\documentclass[11pt]{article}
\usepackage{amssymb,amsfonts,amsthm, amsmath, color, float, mathtools, booktabs, url}
\usepackage{tikz}
\usetikzlibrary{shapes.geometric, arrows}
\usepackage[unicode,psdextra]{hyperref}
\newcommand\pdfmath[1]{\texorpdfstring{$#1$}{#1}}
\usepackage{stackengine,amsmath}
\stackMath
\usepackage{scalerel}
\usepackage{mathrsfs}
\usepackage{caption}
 \def\obrace{\iftrue{\else}\fi}
 \def\cbrace{\iffalse{\else}\fi}
 
 \let\originalparagraph\paragraph
 \renewcommand{\paragraph}[2][.]{\originalparagraph{#2#1}}
 \newcommand{\ol}{\overline}

 \newcommand{\pp}{{\mathbb P}}

 \newcommand{\witi}{\widetilde}
 
 \newcommand{\zz}{{\mathbb Z}}

 \newcommand{\rr}{{\mathbb R}}

 \newcommand{\calc}{{\mathcal C}}

 \newcommand{\calg}{{\mathcal G}}

 \newcommand{\calf}{{\mathcal F}}

 \newcommand{\mz}{{\mathfrak z}}

 \newcommand{\veps}{\varepsilon}

 \newcommand{\beq}{\begin{eqnarray*}}
 	\newcommand{\feq}{\end{eqnarray*}}
 \newcommand{\beqn}{\begin{eqnarray}}
 \newcommand{\feqn}{\end{eqnarray}}
 \newcommand{\bes}{\begin{split}}
 	\newcommand{\fes}{\end{split}}
 \newcommand{\besa}{\begin{align}}
 	\newcommand{\fesa}{\end{align}}
 \newcommand{\as}{\mbox{a.\,s.}}

 \newtheorem{theorem}{Theorem}
 
 \newtheorem*{conj*}{Conjecture}
 \makeatletter \@addtoreset{theorem}{section}\makeatother
 
 \newcommand{\nn}{{\mathbb N}}
 
 \makeatletter \@addtoreset{theorem}{section}\makeatother

 \makeatletter \@addtoreset{theorem}{section}\makeatother
 
 \newtheorem{lemma}[theorem]{Lemma}
 \newtheorem{assume}[theorem]{Assumption}
 \newtheorem*{theorem*}{Theorem}
 
 \newtheorem{proposition}[theorem]{Proposition}
 \newtheorem{corollary}[theorem]{Corollary}
 \newtheorem{remark}[theorem]{Remark}
 \newtheorem{example}[theorem]{Example}
 \makeatletter
 \newcommand{\leqnomode}{\tagsleft@true\let\veqno\@@leqno}
 \newcommand{\reqnomode}{\tagsleft@false\let\veqno\@@eqno}
 \makeatother
 \usepackage{algorithm}
 \usepackage{graphicx}
 \usepackage[noend]{algpseudocode}
 \makeatletter
 \def\BState{\State\hskip-\ALG@thistlm}
 \makeatother

 \DeclareCaptionLabelFormat{noname}{#2}
 \newlength\myindent
 \DeclareMathOperator{\argmax}{\arg\max}

\tikzstyle{startstop} = [rectangle, rounded corners, minimum width=1cm, minimum height=0.5cm,text centered, draw=black, fill=red!30]

\tikzstyle{io} = [trapezium, trapezium left angle=70, trapezium right angle=110, minimum width=1cm, minimum height=0.5cm, text centered, draw=black, fill=blue!30]

\tikzstyle{process} = [rectangle, minimum width=1cm, minimum height=0.5cm, text centered, draw=black, fill=orange!30]
\tikzstyle{decision} = [diamond, minimum width=1cm, minimum height=0.5cm, text centered, draw=black, fill=green!30]

\tikzstyle{arrow} = [thick,->,>=stealth]

 \title{Balancing art and money in pursuit of a Kelly-type optimality}

 \author{
 R.~Rastegar\thanks{Occidental Petroleum Corporation, Houston, TX 77046, USA; e-mail:  reza\_rastegar2@oxy.com}
 	\and
 	A.~Roitershtein\thanks{Department of Statistics, Texas A\&M University, College Station, TX 77843, USA;
 		\newline e-mail: alexander@stat.tamu.edu}
 \and
 V.~Roytershteyn\thanks{Space Science Institute, Boulder, CO 80301, USA;
 		e-mail: vroytershteyn@spacescience.org}
 \and
 	V.~Seetharam\thanks{Department of Computer Science, Texas A\&M University, College Station, TX 77843, USA;
 		\newline e-mail: seetharamg@tamu.edu}}

\begin{document}
 	\maketitle
\begin{abstract}
We introduce and study a mathematical model of an art collector. In our model, the collector is a rational agent whose actions in the art market are driven by two competing long-term objectives, namely sustainable financial health and maintaining the collection. Mathematically, our model is a two-dimensional random linear dynamical system with transformation matrix of a peculiar type. In some examples we are able to show that within the Kelly-type optimization paradigm, that is optimizing the system's Lyapunov exponent over a set of policy parameters, the dilemma ``art or money" can be successfully resolved, namely the optimal policy creates a coexistence equilibrium where the value of both is increasing over the time.

\end{abstract}
 	\noindent{\em MSC2010: } Primary 91B74; 37N40; 37M25  Secondary: 91B70; 37H12; 91B32  \\
 	\noindent{\em Keywords}: art markets, Kelly optimization, Lyapunov exponent, passion investment, random iterations, portfolio optimization
 \section{Introduction}
\label{intro}
\paragraph{Collecting as an economic behavior}
Individual collecting is a routine human behavior which from the economics perspective manifests itself as a blend of luxury consumption with a form of \textit{alternative investment} \cite{belk, investa, route, kleine}. Collecting is a growing phenomenon, its global surge, in particular among Generation~Y, has been facilitated in the past two decades by networking of new information, media and communication technologies \cite{colle, compuls, milad} as well as the emergence of new markets \cite{garay, aus, nocor}. Previous studies indicate that one out of every three Americans might be collecting something \cite{betn, motive}. Credit Suisse estimated in 2020 the value of its client's privately owned collectibles at over US\$1 Trillion in 2020 \cite{ratiop}. In 2018 the size of the art market was estimated to be US\$67.4 billion worldwide \cite{garay}, illuminating its increasing role as a niche market for financial investment, accepted as legitimate by mainstream finance and by economists \cite{coslard, emo5, ratiop, milad, fina}. According to a survey conducted by Barclays in 2012, high-net-worth individuals had in average almost 2\% of their wealth invested in artworks and 8\% in other \textit{emotional assets} or \textit{investments of passion} \cite{coslard, emo5}.
\par
Collecting reflects a multitude of sociological, psychological, economic, and possibly biological motives, including leisure, aesthetics, competition, fantasy, prestige, sensual gratification, extending the self, legacy, and more \cite{belk4, colle, martin, compuls, pear}. Collector's demographics as well as their peculiar personality traits have been examined in \cite{belk, compuls, pear} and more recently in \cite{kleine} and \cite{betn}, the latter reference being specifically focusing on collectors with economic motives. A sociological portrait of art collectors along with a review of the literature on this topic can be found in \cite{compuls, socio}, and a more historical and economic-anthropological perspective on collectors and their motives is offered in \cite{route, pearce}; see also \cite{colle} for a systematic review of various aspects of collecting behavior.
\par
It has been observed by many researchers that ``the distinction between purchasing collectibles for pleasure or for pecuniary reasons is ambiguous and thus difficult to dissociate" \cite{ratiop}. For instance, the authors of \cite{colle} single out Financial Value as one of six main themes in in their classification of collector's motivations (five others are Achievement through Collecting Goals; Social Membership; Cooperation and Competition; Societal and Personal Memories; Legacy) and remark that ``it would be prudent to recognise the financial value of collecting and acknowledge it as a potential motivation for consumers in combination with psychological motivations". A category of collectors that we are interested in are individuals who buy and sell art in ways meant to increase both their status and financial position, they are referred to as \textit{passionate investors} in \cite{isra}, \textit{investor collectors} in \cite{kleine}, \textit{inquisitive collectors} in \cite{compuls}, and  \textit{hierarchically oriented collectors} in \cite{socio}. For example, in a 2013 survey in Germany, initially  fielded with 316,500 panelists, among 4042 responded individuals 225 indicated consistently in their answers that they hold collectibles at least partly for investment purposes \cite{kleine}. When viewed as investors, this group of individuals is not purely financially motivated and thereby considers collectibles as a form of alternative investment, that is a as way to diversify their investment portfolio \cite{betn}.
\paragraph{The model: variables of interest}
We assume a discrete time, $n=0,1,\ldots,$ and at any time point $n$ describe the state of the collector by two numerical variables, a portion of her capital $X_n$ potentially available for an acquisition in the art market and the current value of her art collection $Y_n.$ Rather than to consider the composition of the financial capital and art collection in details, we focus on the evolution of the \textit{macrostate} $(X_n,Y_n)$ and, specifically, the question of the balancing between two conflicting goals, namely long-term maximization of $Y_n$ and minimization of the probability of a \textit{drawdown event} (that is, exhausting to a nearly zero level) of $X_n.$ In the model that we propose, these two goals can in principle be reconciled, in particular the associated art collection process can be sustainingly self-funding in a realistic market behavior scenario.
\par
We assume that while $X_n$ is measured in a regular currency units, say American dollars, $Y_n$ embraces both emotional and financial value of the art collection in the eyes of the owner, and is measured in some different intrinsic units. We are making a simplifying condition
that there is no re-stocking of the initial financial capital $X_0$ and that the art collector is a rational agent in the art market acting to enhance her collection. The variable $X_n$ can be thought of as the balance at time $n$ in a bank account servicing the art collector's operation.
\par
The art market is an example of a market for ``singular goods" whose value in the eyes of the collector is determined by aesthetic judgments rather than by any commonly accessible measurable metric \cite{karpik, emo3}.  In contrast to a financial asset which has a market price determined by a demand and supply mechanism, a unique artwork is valued differently by different potential owners. The magnitude of this \textit{private value} \cite{gaga} is determined by the interaction of the artwork with collector's emotional world, the strength of its attachment to her socio-technological world, potential resale benefits \cite{lovo, quality, ecoa}. The latter, financial part of the valuation, is itself endogenously related to the distribution of tastes and information among potential buyers at the time of resale and therefore can hardly be determined objectively and with certainty \cite{price, class1}.
\par
The enjoyment associated with art ownership is multi-faceted, the emotional benefits that individuals derive from owning an artwork can take different forms, including viewing pleasure, the admiration of artistic skill or genius, feeling connected to the artist, the pride, appropriation of the art object's ``aura",  and social signaling \cite{ecoa, emo3}. Moreover, collector's decisions are often driven by unconscious emotions (such as excitement, anxiety and denial), fantasies, needs, and desires \cite{ecoa}.
\par
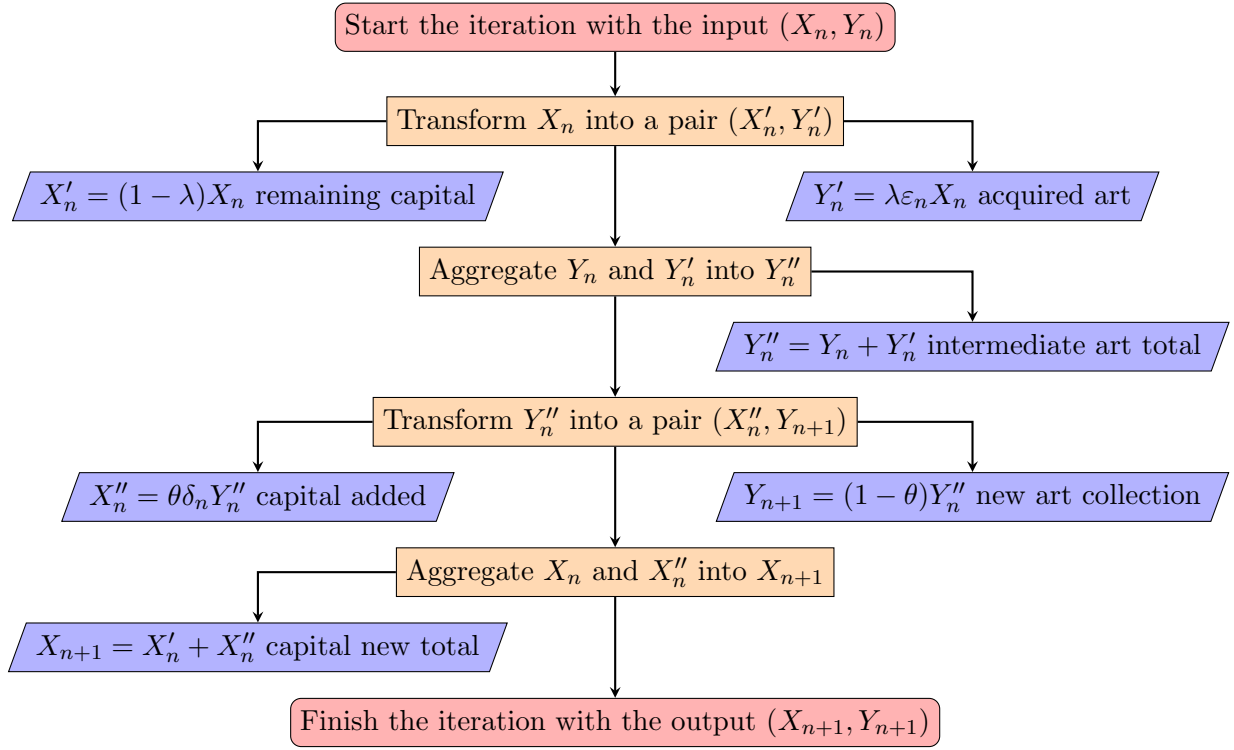
\begin{figure}[!ht]
\centering
\begin{tikzpicture}[node distance=1cm]
\node (start) [startstop] {Start the iteration with the input $(X_n,Y_n)$};
\node (pro1) [process, below of=start,yshift=-0.25cm] {Transform $X_n$ into a pair $(X_n',Y_n')$ };
\node (out1) [io, below of=pro1, xshift=-4.75cm] {$X_n'=(1-\lambda)X_n$ remaining capital};
\node (out2) [io, below of=pro1, xshift=4.75cm] {$Y_n'=\lambda \veps_nX_n$ acquired art};
\node (pro2) [process, below of=out2, xshift=-4.75cm] {Aggregate $Y_n$ and $Y_n'$ into $Y_n''$};
\node (out4) [io, below of=pro2, xshift=4.75cm] {$Y_n''=Y_n+Y_n'$ intermediate art total};
\node (pro3) [process, below of=out4, xshift=-4.75cm] {Transform $Y_n''$ into a pair $(X_n'',Y_{n+1})$};
\node (out5) [io, below of=pro3, xshift=-4.75cm] {$X_n''=\theta \delta_nY_n''$ capital added};
\node (out6) [io, below of=pro3, xshift=4.75cm] {$Y_{n+1}=(1-\theta)Y_n''$ new art collection};
\node (pro4) [process, below of=out6, xshift=-4.75cm] {Aggregate $X_n$ and $X_n''$ into $X_{n+1}$};
\node (out7) [io, below of=pro4, xshift=-4.75cm] {$X_{n+1}=X_n'+X_n''$ capital new total};
\node (stop) [startstop, below of=out7, xshift=4.75cm] {Finish the iteration with the output $(X_{n+1},Y_{n+1})$};
\draw [arrow] (start) -- (pro1);
\draw [arrow] (pro1) -| (out1);
\draw [arrow] (pro1) -| (out2);
\draw [arrow] (pro2) -| (out4);
\draw [arrow] (pro3) -| (out5);
\draw [arrow] (pro3) -| (out6);
\draw [arrow] (pro4) -| (out7);
\draw [arrow] (pro1) -- (pro2);
\draw [arrow] (pro2) -- (pro3);
\draw [arrow] (pro3) -- (pro4);
\draw [arrow] (pro4) -- (stop);
\end{tikzpicture}
\caption{The flowchart of the $(n+1)$-th iteration of the process. By ``transforming" a homogeneous portfolio of assets (either all capital or all art), we mean converting a portion of the portfolio into assets of a different kind while retaining the remaining part. The policy variables $(\lambda,\theta)$ are assumed to be determined at the beginning of the process and remain unchanged during its evolution.}\label{chart}
\end{figure}

The emotions contribute to the part of the utility that sometimes referred in the literature as \textit{emotional dividends} \cite{emo5, lovo, milad} or \textit{psychic returns} \cite{emo5, investa3, frey, frey7}. Without these emotional benefits individuals would be reluctant to allocate their resources to artwork, given other alternatives (such as, for instance, bonds, equities, and real estate) \cite{memoses, geom, samba}. Furthermore, the existence of the non-financial component of the potential buyer's utility has been used to explain certain peculiarities of the art market, such as, for instance, calendar anomalies \cite{investa3, frey7, emo1}, set completion prices \cite{whyc}, and the reluctance of dealers in contemporary art to consider outright price decreases even when it would serve the profit maximization goal. For example, it is argued in \cite{komar, mandel, irpa} that the subordinate role of the price elasticity concept in everyday models of the art market and taboo on decreasing prices can be understood by taking into account, respectively, the observation that collectors infer judgments about the quality of the artwork from its relative price or from a price change and the conspicuous consumption phenomenon.
\paragraph{The model: steps of the stochastic process}
The stochastic process $(X_n,Y_n)_{n\geq 0}$ is assumed to be governed as follows by two \textit{policy parameters} $\lambda \in (0,1)$ and $\theta\in (0,1).$ The parameters are fixed, not changing with time once they are chosen at the beginning of the process. At the beginning of the $n$-th period of time, the collector is investing a portion $\lambda X_n$  of her available wealth $X_n$ in the art market. This investment is translated into the addition of $\lambda X_n \veps_n$ to her collection's value, where $\veps_n$ is a non-negative random variable. To finance her future purchases, the collector is then selling a portion $\theta$ of her art collection which is worth now $Y_n+\lambda X_n \veps_n.$ The profit from this sell, $\theta (Y_n+\lambda X_n\veps_n)\delta_n,$ where $\delta_n$ is another non-negative random variable is then added to the portion of her wealth $(1-\lambda)X_n$ which was put aside at the beginning of the time period; cf.~Fig.\ref{chart}.
\\
\par
Thus, the evolution of the process $(X_n,Y_n)$ is governed by the matrix equation
\beqn
\label{eqa}
\begin{pmatrix}
X_{n+1}
\\
Y_{n+1}
\end{pmatrix}
=
M_n
\begin{pmatrix}
X_n
\\
Y_n
\end{pmatrix}
,
\feqn
where
\beqn
\label{emn}
M_n=\begin{pmatrix}
1-\lambda +\lambda\theta \veps_n \delta_n& \theta\delta_n
\\
\lambda (1-\theta) \veps_n&1-\theta
\end{pmatrix}.
\feqn
The underlying random matrices can be factorized as follows:
\beqn
\label{facts}
M_n=\begin{pmatrix}
\theta \delta_n& 1
\\
1-\theta&0
\end{pmatrix}
\begin{pmatrix}
\lambda \veps_n& 1
\\
1-\lambda&0
\end{pmatrix}.
\feqn
The factorization can be directly related to the above verbal description of the two-step process, see also the flowchart depicted in Fig.~\ref{chart}. Throughout the paper we assume that
\beqn
\label{deha}
\lambda+\theta>0,
\feqn
precluding the trivial case where both $X_n$ and $Y_n$ were never altered.
\paragraph{Dual, $(Y,X)$ versus $(X,Y),$ representation of the model}
The following modified description of the model highlights a symmetry between the $X_n$ and $Y_n$ variables, and in particular shows that the order of the ``buy" and ``sell" operations in Fig.~\ref{chart} is inconsequential for the study of mathematical properties of the model. Suppose that $\theta\neq 1$ and let $\witi X_0=(1-\lambda)X_0,$ $\witi Y_0=Y_0+\lambda \veps_0X_0.$ Denote by $\witi X_n$ and by $\witi Y_n,$ respectively, the value of the wealth and the art value at the beginning of the $n$-th period of time in a variant of our model where the collector is first selling a portion of her collection valued at $\theta \witi Y_n$  to convert it into a financial capital valued $\theta \witi Y_n \delta_n$ and then invest the portion $\lambda$ of the resulting $(\witi X_n+\theta \witi Y_n \delta_n)$-valued financial resource to acquire new art valued $\lambda(\witi X_n+\theta \witi Y_n\delta_n)\veps_{n+1}.$ Thus,
\beq
\witi X_{n+1}&=& (1-\lambda)(\witi X_n+\theta \witi Y_n\delta_n)
\\
\witi Y_{n+1}&=&(1-\theta)\witi Y_n+\lambda (\witi X_n+\theta \witi Y_n\delta_n\veps_{n+1}).
\feq
That is,
\beqn
\label{mtilde}
\begin{pmatrix}
\witi Y_{n+1}
\\
\witi X_{n+1}
\end{pmatrix}
=
\witi M_n
\begin{pmatrix}
\witi Y_n
\\
\witi X_n
\end{pmatrix}
,
\feqn
where
\beq
\witi M_n=\begin{pmatrix}
1-\theta+\lambda \theta \delta_n\veps_{n+1}&\lambda\veps_{n+1}
\\
(1-\lambda)\theta \delta_n&1-\lambda
\end{pmatrix}
=
\begin{pmatrix}
\lambda \veps_{n+1}& 1
\\
1-\lambda&0
\end{pmatrix}
\begin{pmatrix}
\theta \delta_n& 1
\\
1-\theta&0
\end{pmatrix}
.\feq
It is easy to verify that
\beqn
\label{tilde}
\witi X_n=(1-\lambda)X_n
\qquad \mbox{\rm and} \qquad
\witi Y_n= Y_{n+1}(1-\theta)^{-1}.
\feqn
We will use this dual construction in Section~\ref{arm}.
\paragraph{Interpretation of the $(\veps_n,\delta_n)$ variables}
The random variables $\veps_n$ and $\delta_n$ serve as forecasted exchange rates between the idiosyncratic private value and monetary market price for, respectively, buying and selling artworks. The following are among several factors contributing to the intricate variability (in time as well as across the population) and random nature of these elusive indexes. First, art prices are correlated to some degree with the general economy state, they tend to decline during crises, they are also affected by fashions and fads \cite{belk, kleine}. Secondly, collector's own tastes and enthusiasm for collectibles might evolve and fluctuate with time \cite{emo5}. Moreover, bidding behavior and consequently auction prices are impacted by immediate emotions (for instance, community and competition related) and moods of the bidders \cite{auce, tmood}. Thirdly, art market is characterized by heterogeneity, low liquidity, limited and asymmetric (enabling insider's edge) information \cite{hetero, irpa}. Strong evidences supporting the hypothesis of the art markets inefficiency have been found in several studies; see \cite{arte, spectrum} and references therein. Additionally, collector's estimation of an artwork's potential monetary value often relies on quality signals emerging from experts such as gallery owners, curators, art dealers, and critics \cite{price}, which along with endowment effect (overvaluation of an art object owned) and the sunk cost effect (past efforts of building up a collection tend to become a part of its valuation), routinely leads to an inadequate perception of the collection's economic value by the owner \cite{overca, belk4, frey7, vest, amateur}.  Lastly, collectibles are subject to costs and risks such as maintenance, storage, theft, counterfeiting or physical destruction, which apply to a much lesser extent to traditional investment assets \cite{betn, ratiop}.
\par
In view of a low liquidity of art markets \cite{arte, coslard, emo5, komar}, it would be interesting to study a modification of the model with a sparser and perhaps more sporadic over the time market activity. One possibility to achieve such an effect is to make the parameters $(\lambda, \theta)$ into a stochastic process. For instance, at time $n,$ one can replace $\lambda$ with $\lambda i_n$ and $\theta$ with $\theta j_n,$ where $i_n$ and $j_n$ are, possibly correlated, Bernoulli variables. This is the topic of our forthcoming paper \cite{new}.
\paragraph{Dependence structure of the sequence $(\veps_n,\delta_n)_{n\geq 0}$}
In the bulk of the paper we suppose that the pairs $(\veps_n,\delta_n)$ form a stationary ergodic sequence and do not make any further
specific assumptions on its dependence structures (see Assumption~\ref{asu3} below for details). Only in Section~\ref{egp} and examples throughout the paper, we assume that $(\veps_n,\delta_n)_{n\geq 0}$ is a sequence of i.\,i.\,d. pairs (see Assumption~\ref{indi} for details).
Echoing a sentiment expressed in \cite{hamil}, we remark that though a Markov setup is a natural starting point (cf. \cite{geom1, lovo}), nothing in our approach precludes looking at more general probabilistic specifications; see, for instance, \cite{engle, frac1, mendes} for some examples of non-Markovian econometric time series models. For purely technical reasons, it will be convenient to extend $(\veps_n,\delta_n)_{n\geq 0}$ into a double-sided stationary and ergodic sequence $(\veps_n,\delta_n)_{n\in \zz}.$
\paragraph{Optimal policy}
It is observed in \cite{compuls} that while the motivations behind collecting are complex and multifaceted, many of the major motives revolve the development of a more positive sense of self. The authors of \cite{compuls} further postulate that ``collectors are drawn to collecting as a means of bolstering the self by setting up goals that are tangible and attainable and provide the collector with concrete feedback of progress". Furthermore, it is suggested in \cite{frey7} that collectors systematically deviate from the von Neumann-Morgenstern axioms of rational
behavior and, in particular, from subjective expected utility maximization.
\par
We reinterpret these insights by assuming that our collector, at the ``macro-level" that we consider, is seeking to achieve two long-term goals: 1) to maximize the asymptotic growth rate of $Y_n;$ and 2) to maximize the asymptotic growth rate of $X_n.$ Under mild assumptions on the sequence $(\veps_n,\delta_n)_{n\in \zz}$ (see Assumption~\ref{asu3}), these seemingly contradictory goals can be reconciled by optimizing the \textit{Lyapunov exponent}
\beq
\nu(\lambda,\theta):=\lim_{n\to\infty} \frac{1}{n} \ln X_n=\lim_{n\to\infty} \frac{1}{n} \ln Y_n=\lim_{n\to\infty} \frac{1}{n}  E(\ln X_n)=
\lim_{n\to\infty} \frac{1}{n}  E(\ln Y_n),
\feq
where the identities follow from the fundamental results of \cite{fuks}. We thus assume that the goal
of the collector is to identify an \textit{optimal policy} (not necessarily unique) $(\lambda^*,\theta^*)$ such that
\beqn
\label{conti}
\nu(\lambda^*,\theta^*)=\max_{(\lambda,\theta)\in [0,1]^2} \nu(\lambda,\theta).
\feqn
The model offers a combination of basic realistic features along with a certain degree of analytical tractability, we study it both theoretically and numerically. Technically, our model is a two-dimensional twist on the classical Kelly's capital market scheme \cite{kelly, kbook}. In fact, under our assumptions (cf. \eqref{lira} below),
\beq
\nu=\lim_{n\to\infty} E\Big(\ln \frac{X_{n+1}}{X_n}\Big)=\lim_{n\to\infty} E\Big(\ln \frac{Y_{n+1}}{Y_n}\Big).
\feq
That is, similarly to the Kelly's scheme, the collector is optimizing her (in our case, infinite horizon equilibrium rather than the current) logarithmic utility of the immediate investment return. We notice that \eqref{conti} is a one-time infinite-horizon optimization rather than a control problem where the policy may be adjusted at any given time. Thus, the decision maker we consider is ``naive" or ``myopic" in the sense of \cite{feggs, dnow, mayo} and \cite{bgood, kbook} about taking into account the possibility that their future selves might benefit from adjusting the policy due to either personal changes or evolution of the market as the whole. Interestingly enough, the optimization rule \eqref{conti}, as well as many other Kelly-type criteria maximizing an logarithmic utility's growth rate (see, for instance, \cite{cherr, cadane, laur, how, under, response} for a discussion and interpreations), is a natural example of the ergodicity economics optimization paradigm \cite{peters}.
\paragraph{Our contribution}
In two instances (see Section~\ref{key1} and Section~\ref{spec}), using numerical simulations, we show that our model exhibits a two-dimensional analogy of Kelly's effect, namely the optimal policy $(\lambda^*,\theta^*)$ is an interior point of the parameter's domain, rather than a point on its boundary (cf. \cite{kelly}). More specifically, in these two specific examples we show that
\beqn
\label{kef}
\argmax_{(\lambda,\theta)\in [0,1]^2} \nu(\lambda,\theta) \subset (0,1)^2.
\feqn
To show \eqref{kef}, we use a number of theoretical ``piecemeal" insights into the asymptotic behavior of the model, which are collected in
Section~\ref{rgate}, \ref{boundary}, and \ref{arm}. In addition, some illuminating general lower and upper bounds for $\nu(\lambda,\theta)$
for arbitrary values of the parameters $(\lambda,\theta)$ are obtains in Section~\ref{arm}.
\paragraph{Article's organization}
The rest of the paper is structured as follows. In Section~\ref{egp} we obtain the asymptotic of the expectations $E(X_n)$ and $E(Y_n).$
This is later compared with the asymptotic growth rate $\nu(\lambda,\theta),$ formally introduced in Section~\ref{lyapa}. In general, Section~\ref{rgate} is a survey of fundamental properties of $\nu.$ Section~\ref{kellys} is a brief survey of the literature on Kelly-type models. In Section~\ref{ball} we state some basic properties of the model whereas the discussion of more advanced featured is postponed to Section~\ref{arm}, where are main results are presented. In Section~\ref{boundary} we calculate the Lyapunov exponent $\nu(\lambda,\theta)$ on the boundary of the policy domain, that is when at least one of the parameters $\lambda,\theta$ is either zero or one. Sections~\ref{key1} and~\ref{spec} report our numerical studies for some examples when we can reliably claim that \eqref{kef} takes place. Throughout the paper, auxiliary and technical proofs are deferred to Section~\ref{proofs} which plays the role of an appendix.
\section{Rate of the expected growth}
\label{egp}
In this section we are concerned with the asymptotic behavior of the expected values $E(X_n)$ and $E(Y_n).$
The results are stated in a sequence of propositions detailing the long-term behavior of $E(X_n)$ and $E(Y_n),$ with proofs (being a straightforward calculus exercise) either omitted or deferred to Section~\ref{proofs}. While the strategic decision-making of our collector may not be based on the optimization of an affine utility function, the evolution of the mean values is an illuminating feature of the model. A key takeaway from the results presented in this section is the ultimate dependence of the asymptotic behavior of the means on a single stochastic parameter of the model, namely $\gamma$ introduced in \eqref{gamma} below, in addition to the deterministic policy variables $(\lambda,\theta).$
This observation is particularly interesting because $\gamma$ turns out to govern in a similar way the qualitative behavior of the model under the optimal policies introduced in \eqref{conti}, cf. Proposition~\ref{suc}.
\par
Recall \eqref{emn} and let
\beq
\begin{array}{lcll}
\alpha=E(\veps_n)\qquad &\mbox{\rm and}&\qquad \beta=E(\delta_n),&
\\[0.5mm]
U_n=E(X_n)\qquad &\mbox{\rm and}&\qquad V_n=E(Y_n),& \qquad n \geq 0.
\end{array}
\feq
Throughout the paper we will assume the following trivial non-degeneracy condition:
\beqn
\label{abu}
\alpha>0, \qquad \beta>0, \qquad \mbox{\rm and} \qquad U_0+V_0>0.
\feqn
In this section we will combine it with a condition of independence of the underlying variables:
\begin{assume}
\label{indi}
Let \eqref{deha} and \eqref{abu} hold. Assume in addition that
\begin{itemize}
\item [(i)] The pairs of random variables $(\veps_n,\delta_n),$ $n\geq 0,$ are independent of each other.
\item [(ii)] For each fixed $n\geq 0,$ random variables $\veps_n$ and $\delta_n$ are independent of each other.
\end{itemize}
\end{assume}
The assumptions on the dependence structure of the underlying variables will be considerably relaxed when we
will study the almost sure growth rate of the model, cf. Assumption~\ref{asu3}. Some of the results of this section are relevant even under these more general conditions, see the very end of Section~\ref{ball} and in particular \eqref{ergo1} for details.
\par
It follows from \eqref{eqa} that for all $n\geq 0,$
\beqn
\label{eqa1}
\begin{pmatrix}
U_{n+1}
\\
V_{n+1}
\end{pmatrix}
=
M
\begin{pmatrix}
U_n
\\
V_n
\end{pmatrix},
\feqn
where
\beqn
\label{ema}
M:= E(M_n)=\begin{pmatrix}
1-\lambda +\alpha \beta \lambda \theta &\beta\theta
\\
\alpha\lambda(1-\theta)&1-\theta
\end{pmatrix}
.
\feqn
\par
First, for the sake of completeness, we consider extreme cases when either $\lambda\theta =0$ or $\theta=1.$
A common feature of these trivial scenarios is that one of the sequences $U_n$ and $V_n$ saturates
to a constant limit level rather than grows/decays exponentially. More precisely, it is easy to check that \eqref{eqa1} and \eqref{ema} imply the following:
\begin{proposition}
Let Assumption~\ref{indi} hold.
\label{zera}
\item [(i)] If $\theta=0,$ then, for all $n\geq 1,$
\beq
U_n=(1-\lambda)^n U_0 \qquad \mbox{\rm and}\qquad V_n=\alpha\big(1-(1-\lambda)^n\big) U_0 +V_0.
\feq
In particular, $\lim_{n\to\infty} U_n=0$ and $\lim_{n\to\infty} V_n=\alpha U_0+V_0.$
\item [(ii)] If $\lambda=0,$ then, for all $n\geq 1,$
\beq
U_n=U_0+\beta\big(1-(1-\theta)^n\big) V_0 \qquad \mbox{\rm and}\qquad V_n=(1-\theta)^n V_0.
\feq
In particular,  $\lim_{n\to\infty} U_n=U_0+\beta V_0$ and $\lim_{n\to\infty} V_n=0.$
\item [(iii)] If $\theta=1,$ then, for all $n\geq 1,$
\beq
U_n=(1-\lambda+\gamma\lambda)^nU_0+(1-\lambda+\gamma\lambda)^{n-1}\beta V_0 \qquad \mbox{\rm and}\qquad V_n=0,
\feq
where
\beqn
\label{gamma}
\gamma=\alpha \beta.
\feqn
\end{proposition}
\noindent
The proof of the proposition is straightforward, and therefore omitted.
\par
The money-art-money ``exchange rate" $\gamma$ introduced in \eqref{gamma} turns out to be the key parameter governing the asymptotic behavior of the sequences $U_n$  and $V_n$ also in the general case. Let $\mu$ denote the largest eigenvalue of $M,$ that is
\beqn
\nonumber
\mu&=&\frac{1}{2}\Big(2-\lambda -\theta+\gamma \lambda \theta  +\sqrt{ (\theta-\lambda +\lambda\theta\gamma)^2+4\lambda \theta(1-\theta)\gamma}\Big)
\\
&=&
\label{mu}
\frac{1}{2}\Big(2-\lambda -\theta+\gamma \lambda \theta  +\sqrt{ (\lambda +\theta-\gamma\lambda\theta)^2+4(\gamma-1) \lambda \theta}\Big).
\feqn
Notice that $\mu(\lambda,\theta)=\mu(\theta,\lambda)$ and that
\beqn
\label{rez}
\mu=1~\mbox{\rm when}~\lambda\theta=0, \qquad \mu=1-\theta +\theta \gamma ~\,\mbox{\rm if}\,~\lambda=1,
\qquad \mu=1-\lambda +\lambda \gamma  ~\,\mbox{\rm if}\,~\theta=1.
\feqn
We note in passing that in general, if part (i) in Assumption~\ref{indi} is removed, then with $\gamma$ now understood
as $E(\veps_0\delta_0)$ (cf.~\eqref{gamma1} below),
\beq
\mu=\frac{1}{2}\Big(2-\lambda -\theta+\gamma \lambda \theta  +\sqrt{ (\theta-\lambda +\lambda\theta\gamma)^2+4\lambda\theta(1-\theta)\alpha\beta}\Big).
\feq
If $\lambda\theta \neq 0$ and $\theta\neq 1,$ all entries of matrix $M$ are strictly positive, $\mu$ is its Perron-Frobenius eigenvalue, and
therefore, under Assumption~\ref{indi},
\beqn
\label{ker}
\lim_{n\to\infty}\frac{1}{n}\ln U_n=\lim_{n\to\infty}\frac{1}{n}\ln V_n=\ln \mu.
\feqn
The existence of the limits and the identities in \eqref{ker} is a standard linear algebra result, it is an immediate consequence of the fact
that (even if either $U_0=0$ or $V_0=0$) for some constant $c>0,$ $\frac{x}{c}<U_1<cx$ and $\frac{y}{c}<V_1<cy,$ where $(x,y)^T$ is the Perron-Frobenius eigenvector of $M.$ We have:
\begin{proposition}
\label{kthm}
Suppose that Assumption~\ref{indi} holds, $\lambda\theta \neq 0,$ and $\theta\neq 1.$ Then,
$\mu$ is monotone on both $\lambda$ and $\theta.$ More specifically,
\beq
\frac{\partial \mu}{\partial \lambda}>0,\,\frac{\partial \mu}{\partial \theta}>0\quad \mbox{\rm if} \quad \gamma>1 \qquad \mbox{\rm and} \qquad
\frac{\partial \mu}{\partial \lambda}<0,\,\frac{\partial \mu}{\partial \theta}<0 \quad \mbox{\rm if} \quad \gamma<1.
\feq
Furthermore, if $\gamma=1$ then $\mu=1$ for all values of the parameters $\lambda$ and $\theta.$
\end{proposition}
The proof of the proposition is given in Section~\ref{kproof} below. The monotonicity along with \eqref{rez} imply that for all $(\lambda,\theta)\in (0,1)\times (0,1),$
\beqn
\label{muni}
\mu\in (\gamma,1)~\, \mbox{\rm if}~\,\gamma<1; \qquad \mu=1~\, \mbox{\rm if}~\,\gamma=1; \qquad \mbox{\rm and}\quad \mu\in (1,\gamma)~\, \mbox{\rm if}~\,\gamma>1.
\feqn
Notice that \eqref{muni} can be expressed as
\beqn
\label{rez1}
\min\{1,\gamma\}\leq \mu\leq \max\{1,\gamma\},
\feqn
where, under the conditions of the proposition, the equalities are in fact strict unless $\gamma=1.$ Intuitively, the monotonicity reflects the dichotomic nature of the phase transition at $\gamma=1:$  as far as the expected payoff is concerned, either the collector's enterprise is profitable both  financially and emotionally in which case the larger are ``transaction rates" $\lambda$ and $\theta$ the better, or it is not, in which case avoiding the adversary economic environment whatsoever by setting $\lambda=\theta=0$ is the best policy for the collector.
\par
Our next result highlights the asymptotic mean-field dynamics in the case when $\gamma=1.$
\begin{proposition}
\label{zth}
Suppose that Assumption~\ref{indi} holds and $\gamma=1.$ Then, for any initial values $U_0\geq 0,$ and $V_0\geq 0$ and parameter values $\lambda \in [0,1],$ $\theta \in [0,1],$
we have
\beq
\lim_{n\to\infty} U_n= \frac{\theta}{\theta +\lambda(1-\theta)} \big(U_0+\beta V_0\big)\quad \mbox{\rm and}\quad
\lim_{n\to\infty} V_n=\frac{\lambda(1-\theta)}{\theta +\lambda(1-\theta)}\big(\alpha U_0+V_0\big).
\feq
Moreover, $\alpha U_n+V_n$ as well as its $\beta$-multiplier $U_n+\beta V_n$ is a constant sequence (i.\,e., independent of $n$), while
\begin{itemize}
\item [(i)] $U_n$ is strictly increasing and $V_n$ is strictly decreasing if $\lambda <\frac{\beta \theta }{1-\theta} \frac{V_0}{U_0},$
\item [(i)] $U_n$ is strictly decreasing and $V_n$ is strictly increasing if $\lambda >\frac{\beta \theta}{1-\theta} \frac{V_0}{U_0},$
\item [(i)] both $U_n$ and $V_n$ remain constant if $\lambda =\frac{\beta \theta}{1-\theta} \frac{V_0}{U_0},$
\end{itemize}
where $\frac{\beta \theta}{1-\theta}\frac{V_0}{U_0}$ is understood to be infinity if either $U_0=0$ or $\theta=1.$
\end{proposition}
The proof of the proposition is deferred to Section~\ref{zproof}. We conclude this section with the following immediate consequence of \eqref{mu}.
\begin{proposition}
\label{muc}
Let Assumption~\eqref{indi} holds an denote $\xi=\gamma-1.$ Then,
\beq
\mu=1+
\frac{\xi \lambda \theta}{ \lambda +\theta-\lambda\theta}+o(\xi),
\feq
where we used the standard ``little-o" notation, that is $o(\cdot)$ is a function such that $\lim_{\xi\to 0} \frac{o(\xi)}{\xi}=0.$
\end{proposition}
\section{Long-term rate of growth and logarithmic utility}
\label{rgate}
In the remainder of the paper we adopt the following modification of Assumption~\ref{indi}:
\begin{assume}
\label{asu3}
$\mbox{}$
\begin{itemize}
\item [(i)] The sequence of pairs of random variables $(\veps_n,\delta_n),$ $n\geq 0,$ is stationary and ergodic.
\item [(ii)] There exists a constant $C>0$ such that
\beqn
\label{cis}
P\big(C^{-1}<\min\{\veps_n,\,\delta_n\}\leq \max\{\veps_n,\,\delta_n\}<C\big)=1.
\feqn
\item [(iii)] $X_0$ and $Y_0$ are deterministic (non-random) non-negative numbers, not both are zero.
\end{itemize}
\end{assume}
Note that the marginal distributions in the pair $(\veps_n,\delta_n)$ are not required to be independent under this assumption. Consequently, \eqref{gamma} will be substituted in the remainder of the paper with a more general definition
\beqn
\label{gamma1}
\gamma=E(\veps_0\delta_0),
\feqn
which is consistent with \eqref{gamma} in the particular case when the sequences $(\veps_n)_{n\geq 0}$ and $(\delta_n)_{n\geq 0}$ are independent.
\par
The rest of this section is structured s follows. In Section~\ref{lyapa} we formally introduce the main subject of our study, the Lyapunov exponent $\nu(\lambda,\theta)$ associated with the model. Section~\ref{kellys} contains a brief discussion of the Kelly's optimal policy paradigm for investment and betting. Finally, Section~\ref{ball} introduces some basic properties of the Lyapunov exponent pertinent to our model. In particular, Proposition~\ref{conte} shows that $\nu(\lambda,\theta)$ is a continuous functions of its parameters, and hence the $\arg\max$ in \eqref{kef} is well defined (cf. Corollary~\ref{corol}), Proposition~\ref{duet} states an important symmetry property of $\nu,$ and Proposition~\ref{suc} exhibits a phase transition at $\gamma=1$ in the qualitative behavior of the optimal value $\nu(\lambda^*,\theta^*)$ defined in \eqref{conti}.
\subsection{Underlying Lyapunov exponent}
\label{lyapa}
A classical result of \cite{fuks} (see Corollary to Lemma~2 on p.~462 of \cite{fuks}) ensures that if Assumption~\ref{asu3}  holds and $\lambda,\theta\in (0,1),$ then with probability one the following two limits exist and are equal and finite:
\beqn
\label{liml}
\nu=\nu(\lambda,\theta):=\lim_{n\to\infty}\frac{1}{n}\ln X_n=\lim_{n\to\infty}\frac{1}{n}\ln Y_n.
\feqn
Alternatively \cite{fuks},
\beqn
\label{limasol}
\nu=\lim_{n\to\infty}\frac{1}{n}E(\ln  X_n)=\lim_{n\to\infty}\frac{1}{n}E(\ln  Y_n).
\feqn
The common limit in \eqref{liml} and \eqref{limasol} is the (top) Lyapunov exponent of the sequence of random matrices $M_n$ introduced in \eqref{eqa}, that is, with probability one,
\beqn
\label{limas}
\nu=\lim_{n\to\infty} \frac{1}{n} \ln \max\{X_n,Y_n\}=\lim_{n\to\infty}\frac{1}{n}\ln \|M_n\cdots M_1\|=\lim_{n\to\infty}\frac{1}{n}E\big(\ln \|M_n\cdots M_1\|\big),
\feqn
where $\|\cdot\|$ is any matrix norm in dimension $2.$ Part~(ii) of Assumption~\ref{asu3} ensures \eqref{liml}, whereas the existence of the limits and the identities in \eqref{limas} are guaranteed under a weaker assumption (see Theorem~2 in \cite{fuks}) $E\big(\ln (1+\|M_0\|)\big)<\infty.$ Under this condition, the limits in \eqref{liml} would still exists by the Oseledets multiplicative ergodic theorem (see, for instance, Theorem~1.2 in \cite{mogul}), but they would not be necessarily equal. In the absence of \eqref{liml}, maximizing $\nu$ would mean an infinite horizon optimization of the ``most successful among two assets, monetary $X_n$ and passion $Y_n$", compared to optimizing both assets simultaneously when the equality in \eqref{liml} is in place. We remark that a study of the dependence of $\lim_{n\to\infty} \frac{1}{n}\ln \min\{X_n,Y_n\}$ on the parameters $(\lambda,\theta)$ under more general than Assumption~\ref{asu3} conditions appears to be an interesting direction for a further investigation of our model.
\par
Let
\beq
\Omega=\{(x,y)\in\rr^2:0\leq x\leq 1,\,0\leq y\leq 1\}
\feq
and
\beq
\Omega^\circ=\{(x,y)\in\rr^2:0<x<1,\,0<y<1\},
\qquad \partial \Omega=\Omega\slash\Omega^\circ.
\feq
Using numerical simulations, it is shown in Sections~\ref{key1} and~\ref{spec} that in some examples \eqref{kef} holds true, that is
the optimal policy $(\lambda,\theta)$ is an interior point of $\Omega.$ Heuristically, keeping the least profitable asset in the portfolio might mitigate the risk of the ``gambler's ruin" and thus be beneficial in a long-term. This is compared to the following corollary to Proposition~\ref{kthm}:
\beq
\max_{(\lambda,\theta)\in \Omega}\mu(\lambda,\theta)=\max_{(\lambda,\theta)\in \partial \Omega}\mu(\lambda,\theta)=\sup_{(\lambda,\theta)\in \Omega^\circ}\mu(\lambda,\theta),
\feq
where $\mu$ is defined in \eqref{mu}. We conjecture, but were unable to prove it, that under Assumption~\ref{asu3}, $\argmax_{(\lambda,\theta)\in\Omega}\, \nu(\lambda,\theta)\subset \Omega^\circ$ if $\gamma>1,$ $P(\veps_0<1,\delta_0<1)>0,$ and the distribution of $(\veps_0,\delta_0)$ is sufficiently volatile, for instance ``nearly uniform"
on $(-1,a)\times (-1,b)\subset \rr^2$ for some $a,b>2.$
\subsection{Relation to the Kelly capital growth model}
\label{kellys}
The classical Kelly capital growth model can formally be described as the following one-dimensional version of \eqref{eqa}:
\beqn
\label{z3}
Z_{n+1}=(1-\lambda)Z_n+\lambda \veps_n Z_n,
\feqn
where $Z_n$ is the value of an investment portfolio at time $n,$ $\lambda$ is the parameter reflecting the decision of the investor
on what portion of her capital (i.\,e., $\lambda Z_n$) to invest and what portion (i.\,e., $(1-\lambda) Z_n$) to put off each period of time, and  $\veps_n $ is the factor representing the growth of the investment during the $n$-th time period.  Thus, assuming that $\veps_n$ are i.\,i.\,d. variables with a finite mean, $\ln Z_n$ is a classical random walk:
\beq
\ln Z_{n+1}=\ln Z_n+\ln (1-\lambda +\lambda \veps_n).
\feq
The asymptotic speed of the walk is given by the formula
\beqn
\label{vz}
\nu_z(\lambda):=\lim_{n\to\infty} \frac{1}{n}\ln Z_n=E\Big(\ln \frac{Z_{n+1}}{Z_n}\Big)=E\big(\ln (1-\lambda +\lambda \veps_n)\big).
\feqn
Furthermore, provided that $E(\veps_n)>1$ while $P(\veps_0<1)>0,$ the asymptotic speed $\nu_z$ is maximal at a unique point $\lambda^*\in (0,1)$ such that \cite{breiman, kelly}
\beq
E\Big(\frac{\veps_n-1}{1+\lambda^* (\veps_n-1)}\Big)=0.
\feq
In particular, for any $\lambda \neq \lambda^*,$
\beq
\label{lopt}
\lim_{n\to\infty} \frac{Z_n(\lambda^*)}{Z_n(\lambda)}=+\infty, \qquad \mbox{\rm a.\,s.},
\feq
where $Z_n(\lambda)$ is defined in \eqref{z3} and $Z_n(\lambda^*)$ is  its analogue defined by a similar random recurrence with the same sequence $\veps_n$ and the same initial condition $Z_0(\lambda^*)=Z_0(\lambda),$ but using $\lambda^*$ instead of $\lambda$ as the investment policy.
The criterion for evaluating investment strategy $\lambda$ seeking to maximize $\nu_z$ is sometimes referred to as the \emph{Kelly optimal criterion}. By virtue of \eqref{vz}, this policy can be viewed as the maximization of the logarithmic utility $E(\ln \cdot)$ of the immediate capital return $Z_{n+1}/Z_n.$
\par
Kelly's original article \cite{kelly} offers a somewhat surprising interpretation of the optimal investment policy $\lambda^*$ in terms of the information theory. Breiman \cite{breiman1, breiman} (see also \cite{monthly} for an excellent summary account), provided several important characterizations of $\lambda^*.$  For a more recent modifications focusing on alternative interpretations of Kelly's criterion see, for instance, \cite{ brow, median, nopr}.
\par
Heuristically, optimizing the speed of the random walk $\ln Z_n$ yields the fastest (in average, over possible realizations of the sequence $(\veps_n)_{n\geq 0}$) way to achieve preassigned high values of $Z_n$ (see Theorem~1 in \cite{breiman} for a formal statement).
The often criticized downside of such an aggressive policy is the high risk it bears for the investor.
The risk is manifested in the typically high volatility of the time-series trajectory $(Z_n)_{n\geq 0}$ and the consequent, recurrent discrepancy
between the long-term sustainability of the model and its local short-term performance, see \cite{cherr, response} for an illuminating discussion and partial rebuttal of the criticism. For a recent work augmenting Kelly's criterion with constraints based on utilizing finite-horizon performance metrics and measures of risk-aversion, see \cite{risk1,pnas, risk}. We remark that while the volatility of model's trajectories poses even more severe challenges
in higher dimensions (in particular, the Lyapunov exponent $\nu$ is difficult to estimate even numerically \cite{jurga, pcot, pras}), the aggressive investment strategy it entitles is arguably particularly well-suited for a passion ivestment model discussed in our paper \cite{emo5, kbook, sima}.
\par
Kelly's original model has been extended in several ways over the years, with some of the generalizations being directly motivated by applications to economics, decision theory, and financial mathematics. For instance, \eqref{vz} readily extends to stationary and ergodic sequences $\veps_n.$ For an authoritative review of the literature published prior to 2011, see \cite{how}. Thorp's article \cite{sport} offers a superb guide into properties of the original Kelly's model and its applications; see also \cite{ziemba} for a brief and light summary account. For some interesting recent results and literature review, see, for instance, \cite{nike, automa, risk1, compas}.
\subsection{Basic properties of the Lyapunov exponent}
\label{ball}
We begin with the following observation.
\begin{proposition}
\label{conte}
Let Assumption~\ref{asu3} hold. Then, $\nu(\lambda,\theta)$ is continuous on $\Omega.$
\end{proposition}
\begin{corollary}
\label{corol}
Let Assumption~\ref{asu3} hold. Then, there exists
\beqn
\label{karma}
\nu^*:=\max_{(\lambda,\theta)\in \Omega} \nu(\lambda,\theta).
\feqn
\end{corollary}
The proof of the proposition is given in Section~\ref{pconte} below. While the continuity in the interior of the domain $\Omega$ and on the part of the boundary $\partial \Omega$ where $\lambda\theta=0$ follows immediately from general results (see \cite{conti}), the case when either $\lambda=1$ or $\theta=1$ requires an ad-hoc arguments. The argument that we give in Section~\ref{pconte} relies on the explicit formulae
for $\nu(\lambda,\theta)$ on the boundary along with bounds for it within the interior.
\par
Our next result gives necessary and sufficient conditions for $\nu^*$ introduced in \eqref{karma} to be larger than one. The proof of the following proposition is included in Section~\ref{psuc}.
\begin{proposition}
\label{suc}
Let Assumption~\ref{asu3} hold.
\begin{itemize}
\item [(a)] Suppose that $\gamma\leq 1.$
If in addition, $COV(\veps_n,\delta_n)\geq 0,$ then $\nu(\lambda,\theta)\leq 0$ for all $(\lambda,\theta)\in \Omega.$
\item [(b)] Suppose that $\gamma >1.$ Then, $\nu(\lambda,\theta)>0$ within an open neighborhood of any point of $\Omega$ in the form $(\lambda,1)$ for all $\lambda >0$ small enough.
\end{itemize}
\end{proposition}
The following symmetry result follows immediately from \eqref{facts}. Denote
\beq
\Xi=(\veps_n,\delta_n)_{n\in\zz}\qquad \mbox{\rm and}\qquad \Xi_1=(\delta_n,\veps_{n+1})_{n\in\zz}.
\feq
Thus, the sequence $\Xi_1$ is obtained from the original random environment $\Xi$ by first reversing the roles of the margins $\veps_n$ $\delta_n,$ and then shifting the latter sequence one time unit forward. To emphasize the dependence of the Lyapunov exponent on the distribution of $\veps_n$ and $\delta_n$ we will write $\nu_{_\Xi}$ for $\nu.$
\begin{proposition}
\label{duet}
Let Assumption~\ref{asu3} hold. Then, $\nu_{_\Xi}(\lambda,\theta)=\nu_{_{\Xi_1}}(\theta,\lambda)$ for all $(\lambda,\theta)\in \Omega.$
\end{proposition}
The equality $\nu_{_\Xi}(\lambda,\theta)=\nu_{_{\Xi_1}}(\theta,\lambda)$ can be obtained by removing
\beq
\begin{pmatrix}
\theta\delta_n& 1
\\
1-\theta&0
\end{pmatrix}
\qquad
\mbox{\rm and}
\qquad
\begin{pmatrix}
\lambda \veps_0& 1
\\
1-\lambda&0
\end{pmatrix}
\feq
from the product $M_n\cdots M_0$ considered as a product of $2(n+1)$ factors using \eqref{facts}. Alternatively, at least for $\theta\neq 1,$
it can be derived directly from \eqref{mtilde} and \eqref{tilde}.
\par
We conclude this section with a simple classical upper bound for $\nu.$ Let $\mu_n$ be the spectral norm
of matrix $M_n,$ that is the Perron-Frobenius (largest) eigenvalue of $M_n.$ Then, similarly to \eqref{mu},
\beq
\mu_n&=&\frac{1}{2}\Big(2-\lambda -\theta+\gamma_n \lambda \theta  +\sqrt{ (\theta-\lambda+\lambda\theta\gamma_n)^2+4\lambda \theta(1-\theta)\gamma_n }\Big)
\feq
where
\beqn
\label{gamman}
\gamma_n=\veps_n\delta_n.
\feqn
By virtue of \eqref{limas} and the ergodic theorem, we have
\begin{proposition}
\label{ergo}
Let Assumption~\ref{asu3} hold. Then, $\nu \leq E(\ln \mu_0).$
\end{proposition}
In particular, an analogue of \eqref{rez1} for $\mu_n$ implies that
\beqn
\label{ergo1}
\nu \leq E(\ln^+ \gamma_0),
\feqn
where $\ln^+ x$ stands for $\max\{0,\ln x\}.$
\section{Four boundary regimes and their perturbations}
\label{boundary}
In this section we study the model with parameters $(\lambda,\theta)$ chosen on the boundary of the domain $\Omega.$ In particular, we suitably extend the definition of $\nu$ from the interior $\Omega^\circ$ to the boundary $\partial\Omega$ of the domain $\Omega$ in the cases where the Lyapunov exponent of the two-dimensional linear system is not defined by \eqref{liml}, that is the asymptotic growth rates of $X_n$ and $Y_n$ are different. The explicit formulas for $\nu$ on the boundary are used in the proof of the continuity result in Proposition~\ref{conte} (see Section~\ref{pconte} below) and, in addition, facilitate our analysis of the Kelly effect in Sections~\ref{key1} and~\ref{spec} below  (Example \ref{key} in the latter).
\subsection{Speculator regime: \pdfmath{\theta=1}}
\label{specu}
Suppose that $\theta=1.$ By \eqref{emn},  $Y_n=0$ for all $n\geq 1,$ and, for all $n\geq 2,$
\beq
X_{n+1}=(1-\lambda)X_n+\lambda X_n\gamma_n,
\feq
where $\gamma_n$ is introduced in \eqref{gamman}. Therefore,
\beq
\frac{X_n}{X_1}=\prod_{k=0}^{n-1}(1-\lambda+\lambda \gamma_k).
\feq
Hence, by the ergodic theorem,
\beqn
\label{nu3}
\nu(\lambda):=\lim_{n\to\infty} \frac{1}{n}\ln X_n=E\big(\ln (1-\lambda +\lambda \gamma_0)\big), \qquad \mbox{\rm a.\,s.}
\feqn
Notice that by Jensen's inequality, assuming that $\gamma_0$ is not degenerate (non-constant),
\beq
\nu<\ln\big(1-\lambda +\lambda \gamma\big),
\feq
where $\gamma$ is introduced in \eqref{gamma1}. In particular, $\nu<0$ for all $\lambda>0$ if $\gamma\leq 1.$ Furthermore, by the dominated convergence theorem,
\beq
\frac{\partial \nu}{\partial \lambda}=E\Big(\frac{\gamma_0-1}{1-\lambda +\lambda \gamma_0}\Big),
\feq
and hence
\beq
\frac{\partial^2 \nu}{\partial \lambda^2}<0,\qquad \frac{\partial \nu}{\partial \lambda}(0)=\gamma-1,\qquad \frac{\partial \nu}{\partial \lambda}(1)=
1-E\big(\gamma_0^{-1}\big).
\feq
We summarize the above results as follows.
\begin{proposition}
\label{sp1}
Let Assumption~\ref{asu3} hold and assume, in addition, that $\theta=1.$ Then, the Lyapunov exponent $\nu(\lambda,1)=\nu(\lambda)$ is given by \eqref{nu3}. Furthermore, the following holds for the maximum value of $\nu(\lambda):$
\begin{itemize}
\item [(i)] If $\gamma\leq 1,$ then $\nu(0)=0$ and $\nu(\lambda)<0$ for all $\lambda>0.$
\item [(ii)] If $\gamma>1$ and $E(\gamma_0^{-1})>1,$ then there exists a unique $\lambda^*\in (0,1)$ such that
\beq
E\Big(\frac{1}{1-\lambda^*+\lambda^*\gamma_0}\Big)=1,
\feq
in which case
\beq
\nu(0)=0<\nu(\lambda)<\nu (\lambda^*)=E\big(\ln (1-\lambda^* +\lambda^* \gamma_0)\big).
\feq
for all $\lambda \in (0,\lambda^*),$ and $\nu(1)<\nu(\lambda)<\nu(\lambda^*)$ for all $\lambda\in (\lambda^*,1).$
\item [(iii)] If $\gamma>1$ and $E(\gamma_0^{-1})\leq 1,$ then
\beq
\nu(0)=0<\nu(\lambda)<\nu (1)=E\big(\ln \gamma_0\big).
\feq
for all $\lambda \in (0,1).$
\end{itemize}
\end{proposition}
\par
In the rest of the paper, we will occasionally use the following framework for our examples. Let
$\veps, \delta\in (0,1),$ and $\eta>1$ be three given numbers.
\begin{assume}
\label{asu4}
Suppose that Assumption~\ref{indi} holds and, in addition,
\begin{itemize}
\item all variables $\veps_n$ and $\delta_n,$ $n\geq 0,$ are independent of each other,
\item $\veps_n$ is either $\veps$ or $\eta$ with equal probabilities.
\item $\delta_n$ is either $\delta$ or $\eta$ with equal probabilities.
\end{itemize}
\end{assume}
We will refer to the model in Assumption~\ref{asu4} as {\rm BERN}$(\veps,\delta;\eta).$ This class of Bernoulli models is simple, flexible, arguably contains a range of realistic scenarios for an actual art collection, amenable to exact calculations, and hence efficiently illustrates general ideas,
\begin{example}
\label{example}
For a {\rm BERN}$(\veps,\delta;\eta)$ model, we have
\beq
\gamma=\frac{1}{4}(\veps+\eta)(\delta+\eta).
\feq
Furthermore,
\beq
\Delta := E(\gamma_0^{-1})=
\frac{1}{4}\Big(\frac{1}{\veps\delta}+\frac{1}{\veps\eta}
+\frac{1}{\eta\delta}+\frac{1}{\eta^2}\Big).
\feq
For instance, for {\rm BERN}$(0.95,0.95;1.1),$ $\gamma=1.050625>1$ and $\Delta\approx 0.962089<1.$
Therefore, by Proposition~\ref{sp1}, the $\nu$-maximizer for $\theta=1$ is $\lambda=1,$ and
\beq
\max_{\lambda \in [0,1]} \nu(\lambda,1)=\nu(1,1)=E(\ln \gamma_0 )\approx 0.044314,
\feq
while, by Proposition~\ref{kthm}, $\max_{\lambda\in [0,1]}\mu(\lambda,1)=\mu(1,1)=\gamma=1.050625.$
\par
On the contrary, for {\rm BERN}$(0.75,0.75;1.3),$ $\gamma=1.050625>1$ and $\Delta\approx 1.105194>1.$
Therefore, by Proposition~\ref{sp1}, the $\nu$-maximizer for $\theta=1$ is the unique solution to the equation
\beq
\frac{1}{4}\Big(\frac{1}{1-\lambda^*+\lambda^*\cdot 0.75^2}+\frac{2}{1-\lambda^*+\lambda^*\cdot 0.75\cdot 1.3}+\frac{1}{1-\lambda^*+\lambda^*\cdot 1.3^2}\Big)=1,
\feq
which yields $\lambda^*\approx 0.3305.$ Thus, in this case,
\beq
\max_{\lambda \in [0,1]} \nu(\lambda,1)=\nu(\lambda^*,1)=E\big(\ln(1-\lambda^*+\lambda^*\gamma_0)\big)\approx 0.016326,
\feq
while, similarly to the previous example, $\max_{\lambda\in [0,1]}\mu(\lambda,1)=\mu(1,1)=\gamma=1.050625.$
\par
Notice that in this example, as expected, increasing the volatility in the model while keeping $\gamma$ unchanged leads to an inferior
long-term performance of the model.
\end{example}
We will next analyze $\nu(\lambda,\theta)$ in the asymptotic regime $\theta=1-\xi$ as $\xi\to 0.$
Using the notation $\xi=1-\theta,$ one can rewrite \eqref{emn} as
\beq
M_n&=&\begin{pmatrix}
1 -\lambda +(1-\xi)\lambda\veps_n\delta_n &(1-\xi)\delta_n
\\
0& 0
\end{pmatrix}
 +\xi
\begin{pmatrix}
0&0
\\
\lambda\veps_n& 1
\end{pmatrix}
.\feq
Thus, $X_n\geq \prod_{k=0}^n (1 -\lambda +(1-\xi)\lambda\gamma_k),$ and we have
\begin{proposition}
\label{nup}
Let Assumption~\ref{asu3} hold. Then, for any $\lambda,\xi\in (0,1),$
\beqn
\label{nu4}
\nu(\lambda,1-\xi)\geq E\big(\ln (1 -\lambda +(1-\xi)\lambda\gamma_0)\big).
\feqn
\end{proposition}
Note that the bound in \eqref{nu4} is asymptotically exact, that is, in view of \eqref{nu3}, the right-hand side converges to $\nu(\lambda,1)$ when $\xi$ approaches zero.
\subsection{Hoarder regime: \pdfmath{\theta=0}}
\label{hoard}
Assume now that $\theta=0.$ Then, by \eqref{emn},
\beq
X_{n+1}=(1-\lambda)X_n \qquad \mbox{\rm and} \qquad Y_{n+1}=Y_n+\lambda X_n\veps_n.
\feq
In particular, the evolution dynamics of $X_n$ is deterministic, namely $X_n=(1-\lambda)^nX_0.$ For the art value $Y_n$ we therefore have:
\beq
Y_{n+1}=Y_n+\lambda (1-\lambda)^n\veps_nX_0=Y_0+\lambda X_0\sum_{k=0}^n(1-\lambda)^k\veps_k.
\feq
In particular, $Y_n$ converges almost surely, as $n$ approaches infinity, to the limit
\beq
Y_\infty:=Y_0+\lambda X_0\sum_{k=0}^\infty(1-\lambda)^k\veps_k.
\feq
Since $E(Y_\infty)=Y_0+ X_0 E(\veps_0)<\infty,$ the limit is finite with probability one. Accordingly, we append \eqref{liml} with the
following formal definition
\beqn
\label{nu5}
\nu(\lambda,0)=0\qquad \mbox{\rm for all}~\lambda \in [0,1].
\feqn
Notice that, since $\mu_0=1$ when $\theta=0,$ Proposition~\ref{ergo} provides an asymptotically tight upper bound for $\nu(\lambda,\theta)$ when $\theta \to 0.$ We remark that the case when $\veps_0$ and $\delta_0$ are both close to one corresponds to the \textit{weak disorder} regime studied in \cite{perturb} (see, for instance, \cite{texier} and references therein for later developments).
\subsection{Eccentric: \pdfmath{\lambda=1}}
\label{centric}
Suppose that $\lambda=1.$ Then, \eqref{emn} yields
\beq
X_{n+1}=\theta \big(Y_n+X_n\veps_n\big)\delta_n
\qquad {\rm and}\qquad
Y_{n+1}=(1-\theta)(Y_n+X_n\veps_n).
\feq
It follows that $\frac{Y_{n+1}}{X_{n+1}}=\frac{1-\theta}{\theta \delta_n}.$ Therefore,
\beq
\frac{(1-\theta)X_{n+1}}{\theta \delta_n}=Y_{n+1}=(1-\theta)\big(Y_n+\lambda X_n\veps_n)
=(1-\theta)\Big(\frac{(1-\theta)X_n}{\theta \delta_{n-1}}+\lambda X_n\veps_n\Big),
\feq
which implies
\beqn
\label{x4a}
\frac{X_{n+1}}{X_n}=(1-\theta)\frac{\delta_n}{\delta_{n-1}}
 +\theta\veps_n\delta_n.
\feqn
Let
\beqn
\label{dn}
d_n=\frac{\delta_n}{\delta_{n-1}},\qquad n\in\zz.
\feqn
Using this notation, we can rewrite \eqref{x4a} as an explicit formula for $X_{n+1}$ and $Y_{n+1}:$
\beq
X_{n+1}=X_0\prod_{k=0}^n \big((1-\theta)d_n+\theta\gamma_n\big),\qquad Y_{n+1}=\frac{X_{n+1}(1-\theta)}{\theta\delta_n}.
\feq
It follows from \eqref{x4a} that
\beq
\nu(\theta)&=&\lim_{n\to\infty} \frac{1}{n}\ln X_n=\lim_{n\to\infty} \frac{1}{n}\ln Y_n=E\bigg(\ln \Big(\frac{(1-\theta)\delta_n}{\delta_{n-1}}
 +\theta\veps_n\delta_n\Big)\bigg)
 \\
 &=&
 E\big(\ln (1-\theta+\theta\veps_1\delta_0)\big), \qquad \mbox{\rm a.\,s.}
\feq
The result is of course consistent with \eqref{nu3} and Proposition~\ref{duet}. For future reference, we state the following as a formal remark.
Let
\beqn
\label{zetan}
\zeta_n=\veps_n\delta_{n-1} \qquad \mbox{\rm and}\qquad e_n=\frac{\veps_n}{\veps_{n-1}},\qquad n\in\zz.
\feqn
Note that $\gamma_n$ and $\zeta_n$ are identically distribute under Assumption~\ref{indi}, but not in general.
\begin{remark}
\label{sp4}
Using Proposition~\ref{duet}, the analogues of Propositions~\ref{sp1} and~\ref{nup} for $\lambda=1$ can be
obtain by formally exchanging the roles of $\lambda$ and $\theta$ and replacing $\gamma_0$ with $\zeta_0$ in their statements.
\end{remark}
\subsection{Slow liquidation of an inherited collection: \pdfmath{\lambda=0}}
\label{museum}
Assume now that $\lambda=0.$ In this case, the art collection serves to finance the owner.
This can be, for instance, the case for an inherited art collection whose items are sold per re nata by new owners.
In this extreme case, \eqref{emn} yields
\beq
X_{n+1}=X_n+\theta\delta_n Y_n
\qquad
\mbox{\rm and}
\qquad
Y_{n+1}=(1-\theta)Y_n.
\feq
Therefore, $Y_n=(1-\theta)^nY_0$ and
\beq
X_{n+1}=X_n+\theta (1-\theta)^n\delta_nY_0=X_0+\theta Y_0\sum_{k=0}^n(1-\theta)^k\delta_k.
\feq
It follows that, as $n\to\infty,$ $X_n$ converges almost surely to
\beq
X_\infty:=X_0+\theta Y_0\sum_{k=0}^\infty (1-\theta)^k\delta_k.
\feq
Since $E(X_\infty)<\infty,$ the limit is a finite random variable. Accordingly, we append \eqref{liml} with the following formal definition
\beqn
\label{nu7}
\nu(0,\theta)=0\qquad \mbox{\rm for~all}~\theta \in [0,1].
\feqn
Of course, \eqref{nu7} is consistent with \eqref{nu5} and the symmetry result stated in Proposition~\ref{duet}.
\subsection{Numerical example when ``Kelly's effect"  holds true}
\label{key1}
Recall  {\rm BERN}$(\veps,\delta;\eta)$ framework from Assumption~\ref{asu4}.
Fig.~\ref{fignu} below summarizes the results of our numerical study for {\rm BERN}$(0.3,0.2;2)$ described in detail in Section~\ref{ebn} below (see \eqref{cores} in particular).
\par
Note that while $\gamma=1.2515>1$ and, consequently (cf. Proposition~\ref{kthm}),
$\ln \mu>0$ for all $(\lambda,\theta)\in \Omega^\circ,$ the simulations suggest that $\nu(\lambda,\theta)$ is negative in a sizable portion
of the domain $\Omega^\circ$ in a neighborhood of the point $(\lambda,\theta)=(1,1).$ This is consistent with the the result in part (ii) of Proposition~\ref{sp1}, given the fact that $E(\gamma_0^{-1})=5.2708>1$ and $E(\ln \gamma_0)=\frac{1}{2}\ln \frac{6}{100}<0.$
\par
The computation suggests that the maximum of $\nu(\lambda,\theta)$ in $\Omega$ is attained at
\beqn
\label{maxa}
 \lambda^* = 0.265\pm 0.02\qquad \mbox{\rm and}\qquad \theta^*=0.284\pm 0.02,
\feqn
and is equal to
\beqn
\label{ma1}
\nu(\lambda^*,\theta^*) = 0.0199...
\feqn
This to be compared with
\beqn
\label{ma3}
\max_{(\lambda,\theta)\in \partial \Omega}\nu(\lambda,\theta)=\nu(0.132304,1)=0.0160933,
\feqn
computed using part (ii) of Proposition~\ref{sp1} (see Fig.~\ref{figa1} and Fig.~\ref{figa} below). In Section~\ref{ebn} below we estimated
theoretically the difference between the actual value of  $\nu(\lambda^*,\theta^*)$ and the one computed in simulations and shown in \eqref{ma1}.
The outcome is given in \eqref{er11} below, it is considerably smaller than the gap between the estimates in \eqref{ma1} and \eqref{ma3}, effectively proving that \eqref{kef} holds true for {\rm BERN}$(0.3,0.2;2).$
\par
The simulation study reported in Fig.~\ref{fignu} is done using a simple direct method described in Section~\ref{ebn} below, using a grid of $50\times50$ equidistantly covering the interior of rectangle $\Omega.$ The number of iterations $n=n(\lambda,\theta)$ at each point of the grid is determined by the following stopping rule: stop when $\bigl|\frac{\nu_n-\nu_{n-1}}{\nu_n}\bigr|<10^{-3},$ where $\nu_k=\nu_k(\lambda,\theta)$ is the approximation of the value of $\nu(\lambda,\theta)$ computed at iteration $k.$
\par
We then ran separately 20 iterations of the same method at the point $(\lambda^*,\theta^*)$ given in \eqref{maxa}. The values $\nu_k(\lambda^*,\theta^*)$ for $k=11,\ldots,20$ are as follows:
{\small
\beqn
\nonumber
&&
\nu_{11}=0.019915518445, \quad \nu_{12}=0.019917137593,  \quad \nu_{13}=0.019918021008, \quad
\nu_{14}=0.019918502880,
\\
\nonumber
&&
\nu_{15}=0.019918765678, \quad \nu_{16}=0.019918908984,
\quad \nu_{17}=0.019918987124, \quad \nu_{18}=0.019919029729,
\\
\label{figsa}
&&
\nu_{19}=0.019919052953,\quad
\nu_{20}=0.019919065598.
\feqn
}
The graphs of $\nu_n(\lambda^*,\theta^*)-\nu_{20}(\lambda^*,\theta^*)$ and $\nu_n(\lambda^*,\theta^*)-\nu_{n-1}(\lambda^*,\theta^*)$ (on a logarithmic scale) for this simulation run ar plotted in, respectively, Fig~\ref{figs} and Fig~\ref{fig33} below.
\par
Finally, to further verify our simulations result for {\rm BERN}$(0.3,0.2;2)$  at the point $(\lambda^*,\theta^*)$ given in \eqref{maxa}, we also ran 16 rounds of iterations of the method for computation Lyapunov exponents introduced in \cite{jurga}. The values $\nu_k(\lambda^*,\theta^*)$ for $k=10,\ldots,16$ are as follows:
{\small
\beqn
\nonumber
&&
\nu_{10}=0.020386826647, \quad \nu_{11}=0.019882803069,  \quad \nu_{12}=0.019921362049, \quad
\nu_{13}=0.019918963924,
\\
\label{ju1}
&&
\nu_{14}=0.019919085714, \quad \nu_{15}=0.019919080635,
\quad \nu_{16}=0.019919081020.
\feqn
}
The graphs of $\nu_n(\lambda^*,\theta^*)-\nu_{20}(\lambda^*,\theta^*)$ and $\nu_n(\lambda^*,\theta^*)-\nu_{n-1}(\lambda^*,\theta^*)$ (on a logarithmic scale) for this simulation run ar plotted in, respectively, Fig~\ref{figs1} and Fig~\ref{fig333} below.
\begin{figure}[!ht]
\centering
\includegraphics[scale=0.65]{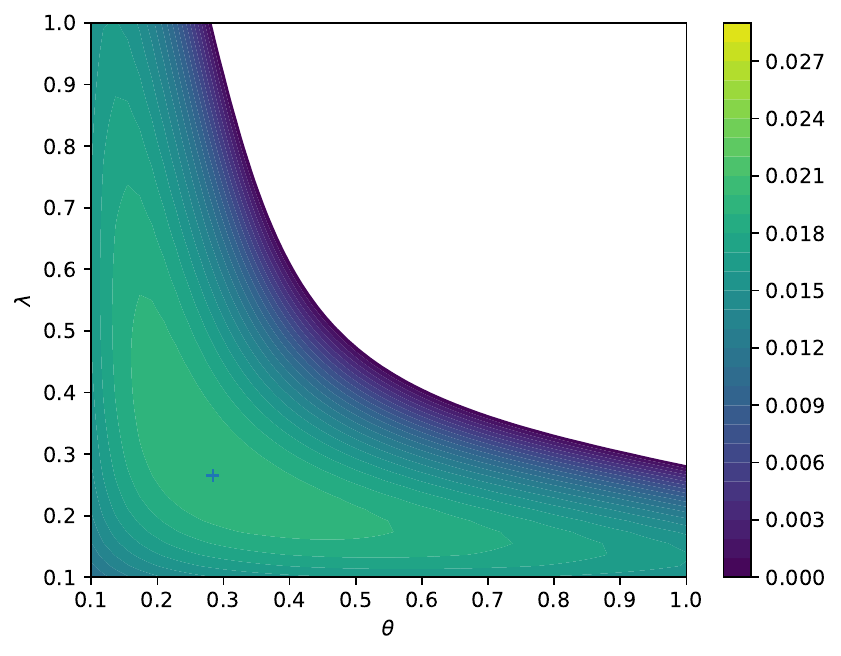}
\caption{The heatmap represents the values of the Lyapunov exponent $\nu$ for the {\rm BERN}$(0.3,0.2;2)$ model for a grid of parameters $(\lambda,\theta)$
equidistantly covering the rectangle $\Omega.$ The white area in the picture corresponds to negative values of $\nu.$ \label{fignu}}
\end{figure}
\begin{figure}[!ht]
\begin{minipage}[b]{0.48\linewidth}
\centering
\includegraphics[width=\textwidth]{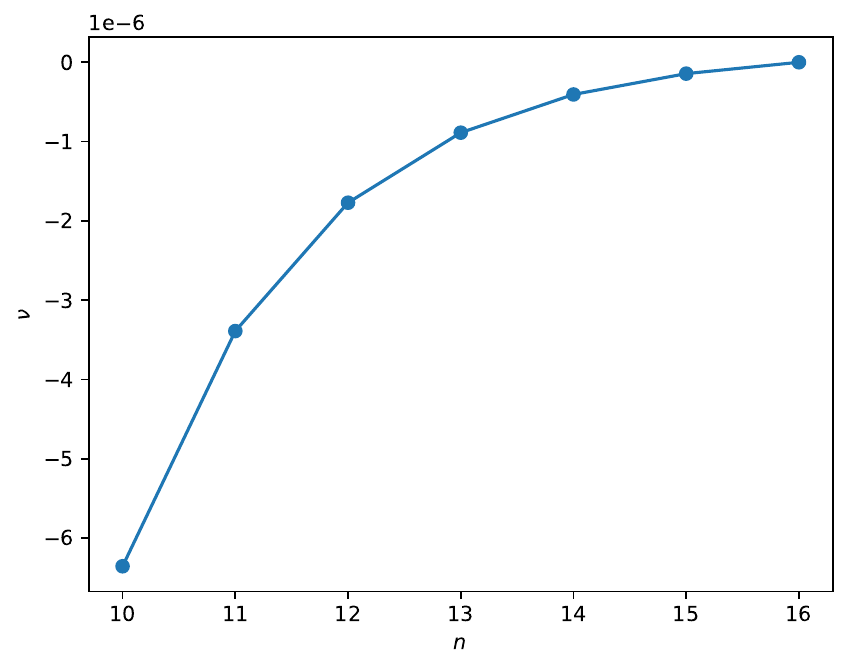}
\caption{A plot of time series differences $\nu_n(\lambda^*,\theta^*)-\nu_{20}(\lambda^*,\theta^*)$  obtained in the simulations reported in \eqref{figsa}.}
\label{figs}
\end{minipage}
\quad
\begin{minipage}[b]{0.48\linewidth}
\centering
\includegraphics[width=\textwidth]{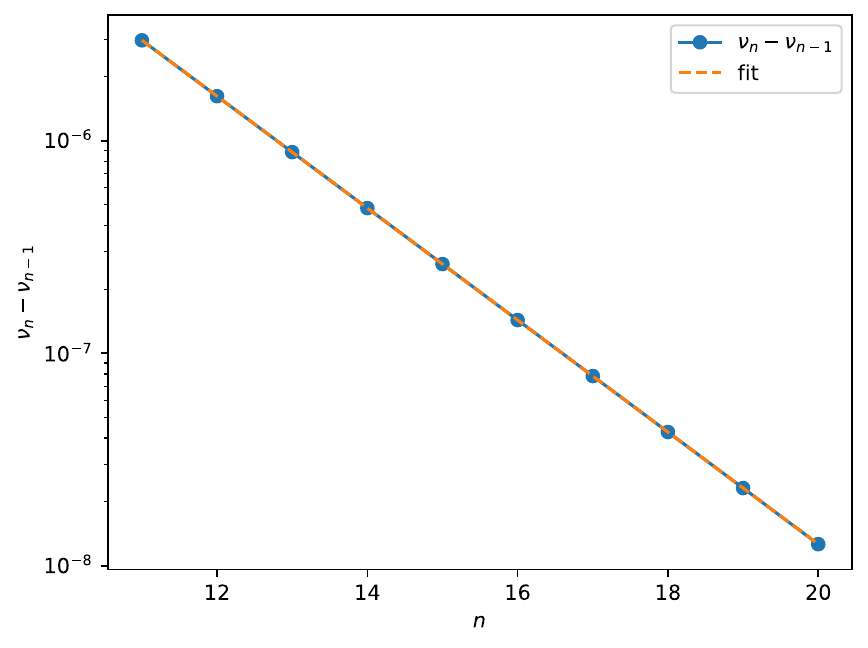}
\caption{Assuming that $\nu_n=\nu_0-\kappa_1e^{-\kappa_2 n},$ we fitted the differences $\nu_n-\nu_{n-1},$ $n=11,\ldots, 20,$ using the least square method, and obtained the estimates $\kappa_2 = 0.606395,$ $\kappa_1=0.00280925.$}
\label{fig33}
\end{minipage}
\end{figure}
\begin{figure}[!ht]
\begin{minipage}[b]{0.48\linewidth}
\centering
\includegraphics[width=\textwidth]{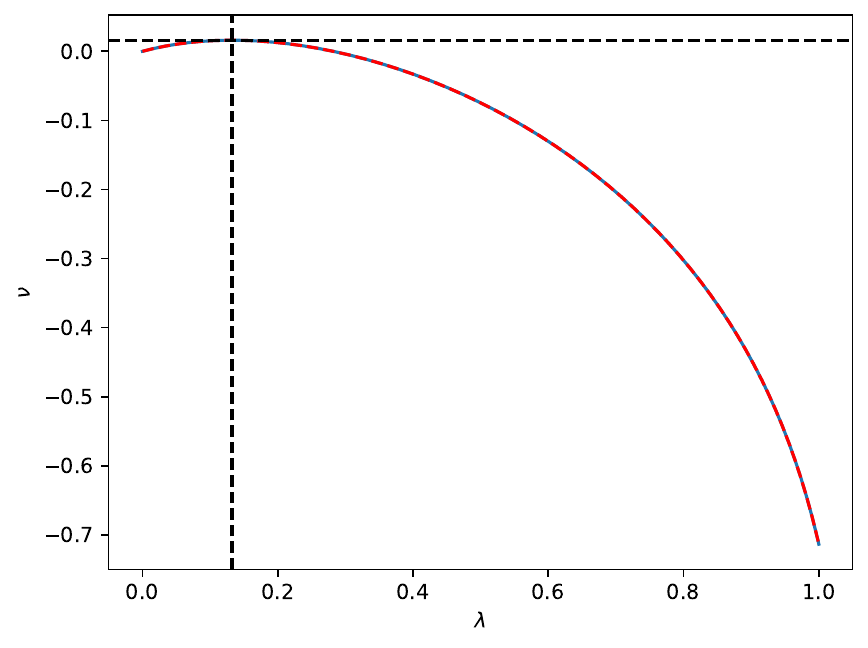}
\caption{Graph of the Lyapunov exponent $\nu(\lambda,1)$ for the {\rm BERN}$(0.3,0.2;2)$ model.}
\label{figa1}
\end{minipage}
\hspace{0.5cm}
\begin{minipage}[b]{0.48\linewidth}
\centering
\includegraphics[width=\textwidth]{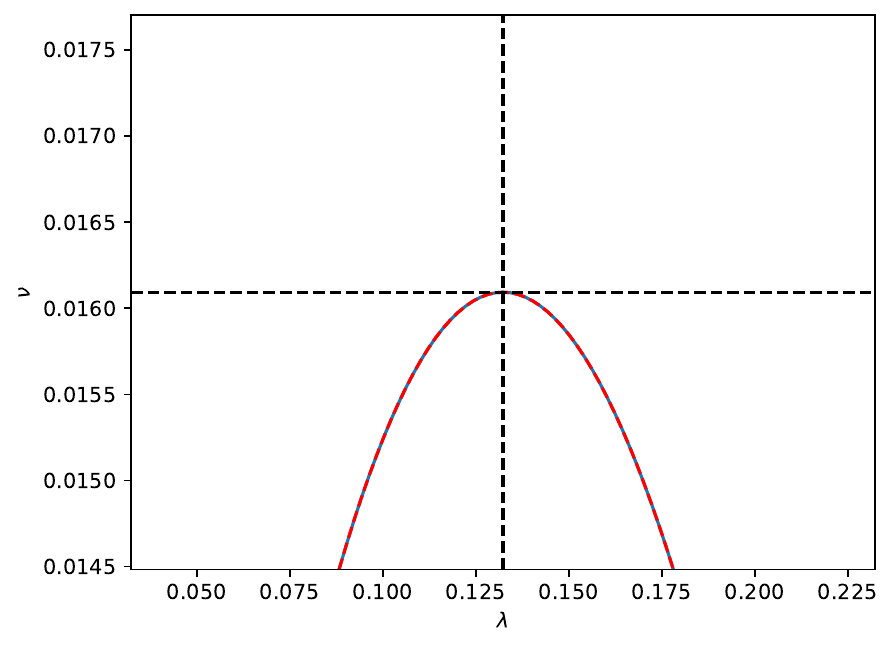}
\caption{Graph of the Lyapunov exponent $\nu(\lambda,1)$ for the {\rm BERN}$(0.3,0.2;2)$ model zoomed-in in a neighborhood of the maximum.}
\label{figa}
\end{minipage}
\end{figure}
\begin{figure}[!ht]
\begin{minipage}[b]{0.48\linewidth}
\centering
\includegraphics[width=\textwidth]{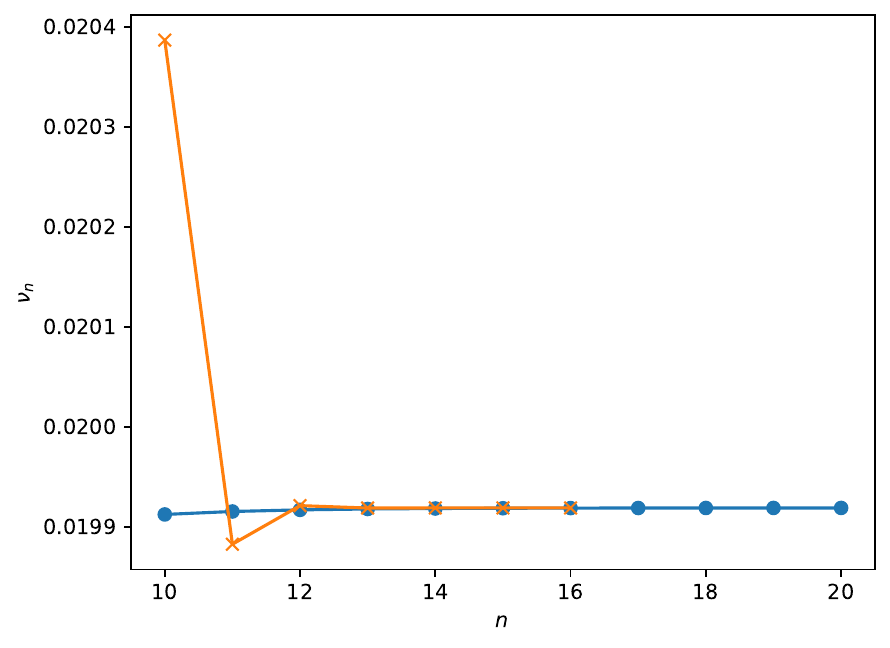}
\caption{Comparison of the simulations results shown in \eqref{figsa} and \eqref{ju1} (the latter is using the method introduced in \cite{jurga}). The yellow line corresponds to the time series
in \eqref{ju1} and the blue line to that in \eqref{figsa}.}
\label{figs1}
\end{minipage}
\hspace{0.5cm}
\begin{minipage}[b]{0.48\linewidth}
\centering
\includegraphics[width=\textwidth]{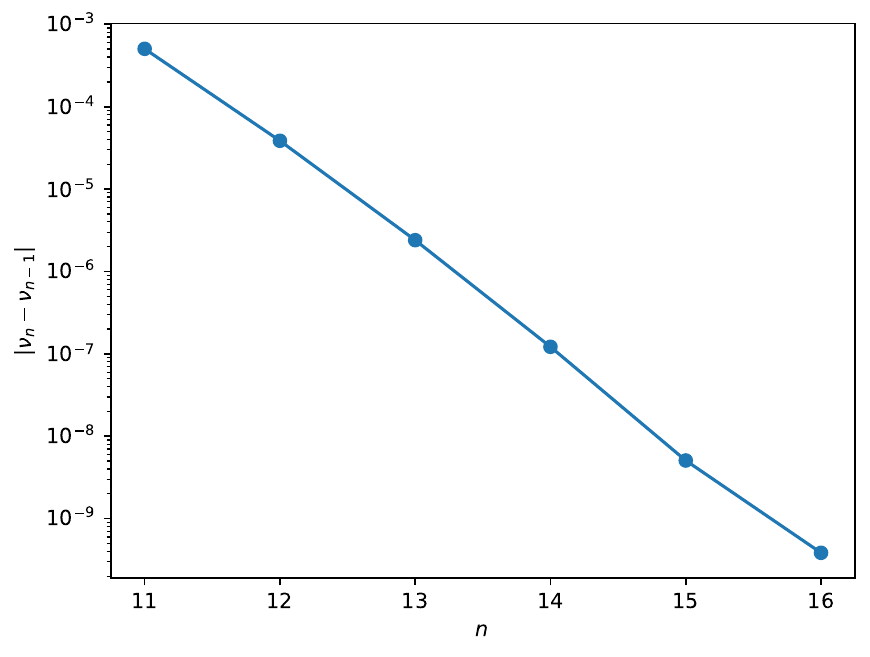}
\caption{The differences $\nu_n-\nu_{n-1},$ $n=11,$ $12,$ $\ldots,$ $16,$ for the simulation results reported in \eqref{ju1}, computing $\nu(\lambda^*,\theta^*)$ using the method introduced in \cite{jurga}.}
\label{fig333}
\end{minipage}
\end{figure}
\section{Upper and lower bounds for \pdfmath{\nu(\lambda,\theta)}}
\label{arm}
The goal of this section is to obtain several bounds for the Lyapunov exponent $\nu.$ These bounds are obtained by using several different techniques, each utilizing an amalgam of a classical approach to the study of products of i.\,i.\,d. matrices with ad-hoc estimates based on the peculiar structure of our matrices.
\par
The distinct structure of the underlying matrices $M_n$ not only enables simultaneous implementation of these classical methods in a single setting, but also allows us to investigate the model in a general ergodic setting of Assumption~\ref{asu3} rather than making a standard i.,i.,d. or Markov-dependence assumptions for the underlying sequence of matrices $(M_n)_{n\geq 0}.$
\par
The rest of this section is organized as follows. In Section~\ref{decouple} we obtain an equivalent representation of the model as a system of two linear second order recursions and exploit it to obtain upper and lower bounds for $\nu(\lambda,\theta)$ (Theorem~\ref{art} and Corollary~\ref{artc}). In Section~\ref{cofr} we use the above recursions to obtain a continued fraction representation for $\frac{X_{n+1}}{X_n}$ $\frac{Y_{n+1}}{Y_n},$ and consequently for $\nu.$ The basic limit results for these fractions is stated in Theorem~\ref{verba}, it yields an exponentially fast converging series representation and new upper and lower bounds for $\nu$ (see Proposition~\ref{ester}, Corollary~\ref{core}, and Proposition~\ref{ester1}). A different method, utilizing a representation of matrices $M_n$ as actions in the projective space and the interplay between continued fractions and linear fractional maps (cd. \cite{beardon}), is employed in Section~\ref{fur} to obtain an alternative continued fraction representation for $\nu,$ the result is stated as Proposition~\ref{nep}. Finally, in Section~\ref{spec} we consider a special example, a modification of the one given in \cite{letac}, where an integral formula (easily computable using standard numerical methods) can be given for the ``diagonal elements" $\nu(\lambda,\lambda).$ Does it give us an example when \eqref{kef} holds true?
\subsection{Decoupling the 2-dim system into 2 autoregressive recursions}
\label{decouple}
The aim of this subsection is to obtain an autoregressive representation \cite{garh} of the model which decouples the underlying two-dimensional dynamical system into two independent second-order linear recursions, see equations \eqref{ar} and \eqref{ar9} below. We note that such a representation is not available in general for two-dimensional linear systems, and its existence for our model reflects a peculiar structure of matrices $M_n.$ The decoupling yields a new continued fraction representation of the Lyapunov exponent $\nu$ (Theorem~\ref{verba}) as well as new explicit bounds for it (Theorem~\ref{art} and Corollary~\ref{artc}).
\par
Write \eqref{emn} explicitly as a system of linear equations:
\beqn
\label{fax}
X_{n+1}&=&(1-\lambda)X_n+\lambda\theta\veps_n\delta_n X_n+\theta \delta_nY_n,
\\
\label{fay}
Y_{n+1}&=&(1-\theta)Y_n+\lambda(1-\theta)\veps_nX_n.
\feqn
Iterating \eqref{fay}, we get
\beq
Y_{n+1}&=&(1-\theta)Y_n+\lambda(1-\theta)\veps_nX_n
\\
&=&
(1-\theta)^2Y_{n-1}+\lambda(1-\theta)^2\veps_{n-1}X_{n-1}+\lambda(1-\theta)\veps_nX_n
\\
&=&\cdots=(1-\theta)^{n+1}Y_0+\sum_{k=0}^n \lambda(1-\theta)^{n+1-k}\veps_kX_k.
\feq
The last identity reflects the fact that each item in the art collection either was a part of the original collection at time zero or
its origin can be traced to a time $k\geq 0$ when it was acquired with funds $\lambda X_k$ and never sold thereafter (an hence the factor $(1-\theta)^{n+1-k}$). Plugging-in the corresponding expression for $Y_n$ into \eqref{fax}, we obtain
\beqn
\label{first}
X_{n+1}&=&(1-\lambda)X_n+\lambda\theta\delta_n\veps_nX_n+\theta \delta_nY_n
\\
\nonumber
&=&
(1-\lambda)X_n +\theta(1-\theta)^n\delta_nY_0+\sum_{k=0}^n \lambda\theta(1-\theta)^{n-k}\delta_n\veps_kX_k
\feqn
Therefore,
\beq
\frac{X_{n+1}-(1-\lambda)X_n}{(1-\theta)^n\delta_n}&=&
\theta Y_0+\sum_{k=0}^n \lambda\theta(1-\theta)^{-k}\veps_kX_k
\\
&=&
\lambda\theta(1-\theta)^{-n}\veps_nX_n+
\frac{X_n-(1-\lambda)X_{n-1}}{(1-\theta)^{n-1}\delta_{n-1}},
\feq
which, using the notation introduced in \eqref{gamman} and \eqref{dn},
can be rewritten as the following linear recursion (\emph{generalized autoregressive model of order $2$ with random coefficients} \cite{garh}):
\beqn
\label{ar}
X_{n+1}=\big(1-\lambda+\lambda\theta\gamma_n+(1-\theta)d_n\big)X_n
-
(1-\lambda)(1-\theta)d_nX_{n-1}.
\feqn
It follows from this recursion that
\beqn
\label{ara}
X_n\leq \big(1-\lambda+\lambda\theta\gamma_{n-1}+(1-\theta)d_{n-1}\big)X_{n-1},
\feqn
and hence
\beq
\frac{X_{n+1}}{X_n}&\leq& 1-\lambda+\lambda\theta\gamma_n+(1-\theta)d_n
-
(1-\lambda)(1-\theta)d_n\big(1-\lambda+\lambda\theta\gamma_{n-1}+(1-\theta)d_{n-1}\big)^{-1}
\\
&=& 1-\lambda+\lambda\theta\gamma_n
+
(1-\theta)d_n\frac{\lambda\theta\zeta_{n-1}+1-\theta} {(1-\lambda)d_{n-1}^{-1}+\lambda\theta\zeta_{n-1}+1-\theta}.
\feq
Consequently,
\beqn
\label{nam1}
\nu \leq E\Big(\ln \Big( 1-\lambda+\lambda\theta\gamma_n
+
(1-\theta)d_n\,\frac{\lambda\theta\zeta_{n-1}+1-\theta} {(1-\lambda)d_{n-1}^{-1}+\lambda\theta\zeta_{n-1}+1-\theta}\Big)\Big).
\feqn
An alternative bound for $\nu(\lambda,\theta)$ can be obtained from \eqref{ar} and \eqref{ara} as follows.
By plugging-in the bound for $X_n$ from \eqref{ara} into the first term in the right-hand side of \eqref{ar}, we obtain
\beq
&& \frac{X_{n+1}}{X_{n-1}}\leq (1-\lambda+\lambda\theta\gamma_n+(1-\theta)d_n\big)
(1-\lambda+\lambda\theta\gamma_{n-1}+(1-\theta)d_{n-1}\big)
-
(1-\lambda)(1-\theta)d_n,
\feq
and, consequently,
{\small
\beqn
\label{nam3}
 \nu^2 \leq E\big(\ln\big\{(1-\lambda+\lambda\theta\gamma_n+(1-\theta)d_n\big)
(1-\lambda+\lambda\theta\gamma_{n-1}+(1-\theta)d_{n-1}\big)
-
(1-\lambda)(1-\theta)d_n \big\}\big).
\feqn
}
On the other hand, it follows from \eqref{first} that
\beq
\frac{X_{n+1}}{X_n}&=&1-\lambda +\lambda\theta\veps_n\delta_n+\theta \delta_n\frac{Y_n}{X_n}.
\feq
Furthermore, by virtue of \eqref{emn},
\beq
\frac{Y_n}{X_n}
&\geq&
\min \Big\{\frac{\lambda (1-\theta)\veps_{n-1}}{1-\lambda +\lambda\theta\veps_{n-1}\delta_{n-1}},\,
\frac{1-\theta}{\theta\delta_{n-1}}\Big\}
\\
&=&\frac{\lambda (1-\theta) \veps_{n-1}}{1-\lambda +\lambda\theta \veps_{n-1}\delta_{n-1}},
\feq
and hence
\beq
\frac{X_{n+1}}{X_n}&\geq &1-\lambda +\lambda\theta\delta_n\veps_n
+\frac{(1-\theta)\lambda \theta \delta_n\veps_{n-1}}{1-\lambda +\lambda\theta \veps_{n-1}\delta_{n-1}}.
\feq
The results of this subsection so far can be summarized as follows:
\begin{theorem}
\label{art}
Let Assumption~\ref{asu3}hold. Then, for any interior point $(\lambda,\theta)\in \Omega^\circ:$
\item [(i)] The sequence $X_n$ obeys the linear recursion \eqref{ar}.
\item [(ii)] The following lower bound holds for the Lyapunov exponent $\nu:$
\beq
\nu \geq E\Big(\ln \Big(1-\lambda +\lambda\theta\gamma_n+(1-\theta)d_n\frac{\lambda \theta \gamma_{n-1}}{1-\lambda
+\lambda\theta \gamma_{n-1}}\Big)\Big).
    \feq
\item [(iii)] The upper bounds in \eqref{nam1} and \eqref{nam3} hold for the Lyapunov exponent $\nu.$
\end{theorem}
Utilizing the notation introduced in \eqref{zetan} and Proposition~\ref{duet}, we obtain
\begin{corollary}
\label{artc}
Let Assumption~\ref{asu3}hold. Then, for any interior point $(\lambda,\theta)\in \Omega^\circ:$
\item [(i)] The sequence $Y_n$ obeys the linear recursion
\beqn
\label{ar9}
Y_{n+1}&=&\big(1-\theta+\lambda\theta\zeta_n+(1-\lambda)e_n\big)Y_n-(1-\lambda)(1-\theta)e_nY_{n-1}.
\feqn
\item [(ii)] The following lower bound holds for the Lyapunov exponent $\nu:$
\beq
\nu \geq E\Big(\ln \Big(1-\theta +\lambda\theta\zeta_n+(1-\lambda)e_n\frac{\lambda \theta \zeta_{n-1}}{1-\theta
+\lambda\theta \zeta_{n-1}}\Big)\Big).
    \feq
\item [(iii)] The following two upper bounds hold for the Lyapunov exponent $\nu:$
\beq
\nu \leq E\Big(\ln \Big( 1-\theta+\lambda\theta\zeta_n
+
(1-\lambda)e_n\,\frac{1-\lambda+\lambda\theta\gamma_{n-1}} {1-\lambda+\lambda\theta\gamma_{n-1}+(1-\theta)e_{n-1}^{-1}}\Big)\Big),
\feq
{\small
\beq
 \nu^2 \leq E\big(\ln\big\{(1-\theta+\lambda\theta\zeta_n+(1-\lambda)e_n\big)
(1-\theta+\lambda\theta\zeta_{n-1}+(1-\lambda)e_{n-1}\big)
-
(1-\lambda)(1-\theta)e_n \big\}\big).
\feq
}
\end{corollary}
\subsection{Continued fraction representation of \pdfmath{\nu}}
\label{cofr}
Recall \eqref{gamman}, \eqref{dn}, and \eqref{zetan}. For $n\in\zz,$ using the double-infinite version of the stationary sequence $(\veps_n,\delta_n),$ let
\beqn
\label{pq}
p_n=1-\lambda+\lambda\theta\gamma_n+(1-\theta)d_n \qquad &\mbox{\rm and}&\qquad
q_n=-
(1-\lambda)(1-\theta)d_n,
\\
\nonumber
r_n=1-\theta+\lambda\theta\zeta_n+(1-\lambda)e_n\qquad &\mbox{\rm and}&\qquad
s_n=-(1-\lambda)(1-\theta)e_n.
\feqn
Thus, for $n\geq 1,$
\beqn
\label{xyrec}
X_{n+1}=p_nX_n+q_nX_{n-1}\qquad \mbox{\rm and}\qquad Y_{n+1}=r_nY_n+s_nY_{n-1}.
\feqn
Equivalently,
\beqn
\label{mxy}
\begin{pmatrix}
X_{n+1}
\\
X_n
\end{pmatrix}
=
\begin{pmatrix}
p_n&q_n
\\
1&0
\end{pmatrix}
\begin{pmatrix}
X_n
\\
X_{n-1}
\end{pmatrix}
\quad
\mbox{\rm and}
\quad
\begin{pmatrix}
Y_{n+1}
\\
Y_n
\end{pmatrix}
=
\begin{pmatrix}
r_n&s_n
\\
1&0
\end{pmatrix}
\begin{pmatrix}
Y_n
\\
Y_{n-1}
\end{pmatrix}
\feqn
Fibonacci-like random recursions in the form similar to \eqref{xyrec} have been studied by many authors, see for instance  \cite{darm} and references therein. The reduction from \eqref{eqa} to \eqref{mxy} enables us to exploit a well-known relation between random Fibonacci sequences and continuous fractions in order to obtain a new representation for the Lyapunov exponent $\nu$ see (Theorem~\ref{verba} below).
\par
For any integer $n\geq 0$ and real $\mz\geq 0$ define continued fractions
{\small
\beqn
\label{ubc}
u_n(\mz)=p_0 + \frac{q_0\kern7em}{\displaystyle
  p_{-1} + \frac{q_{-1}\kern6em}{\displaystyle
    p_{-2} +\stackunder{}{\ddots\stackunder{}{\displaystyle
      {}+ \frac{q_{n-2}}{\displaystyle
        p_{-(n-1)} + \frac{q_{-(n-1)}}{p_{-n}+q_{-n}\mz}}}
}}}
\feqn
and
\beq
v_n(\mz)=r_0 + \frac{s_0\kern5em}{\displaystyle
  r_{-}1 + \frac{s_{-1}\kern4em}{\displaystyle
    r_{-2} +\stackunder{}{\ddots\stackunder{}{\displaystyle
      {}+ \frac{s_{-(n-2)}}{\displaystyle
        r_{-(n-1)} + \frac{s_{-(n-1)}}{r_{-n}+s_{-n}\mz}}}
}}}
\feq
}
Notice that
\beq
p_n+q_n\mz> 0~\mbox{\rm for}~0\leq \mz <(1-\lambda)^{-1}\qquad \mbox{\rm and}\qquad r_n+s_n\mz> 0~\mbox{\rm for}~0\leq \mz <(1-\theta)^{-1}.
\feq
It is easy to verify that both $u_n(\mz)$ and $v_n(\mz)$ are monotone decreasing sequence for, respectively,
$\mbox{\rm for}~0\leq \mz <(1-\lambda)^{-1},$ and $\mbox{\rm for}~0\leq \mz <(1-\theta)^{-1},$ and hence converging
almost surely to a limit within these ranges of the initial value $\mz.$ In fact, once the existence of the Lyapunov exponent $\nu$ has been established independently, the bulk of the argument in \cite{hope} can be carried over to obtain the following result (see also Lemma~\ref{station} below for an extension). .
\begin{theorem}
\label{verba}
Let Assumption~\ref{asu3} hold, and let
\beq
u(\mz)=\lim_{n\to\infty} u_n(\mz)~\mbox{\rm for}~0\leq \mz <(1-\lambda)^{-1} ,\qquad
v(\mz)=\lim_{n\to\infty} v_n(\mz)~\mbox{\rm for}~0\leq \mz <(1-\theta)^{-1}.
\feq
Denote $u_n=u_n(0),$ $v_n=v_n(0),$ $u=u(0),$ and $v=v(0).$ Then,
\beq
\nu=\lim_{n\to\infty} E(\ln u_n)=\lim_{n\to\infty} E(\ln v_n)=E(\ln u)=E(\ln v).
\feq
\end{theorem}
A short proof of Theorem~\ref{verba} is included in Section~\ref{pverba}. The result gives yet another equivalent expression for $\nu.$ The expected values of $\ln u$ and $\ln v$ can be calculated numerically, at least in principle. In particular, though we do not pursue the line of model's investigation in this paper, it seems plausible that
the algorithm introduced in \cite{wright} can be adopted for a class of continuous distributions $(\veps_0,\delta_0).$ In fact, we have the following:
\begin{proposition}
\label{ester}
Let Assumption~\ref{asu3} hold and denote
\beqn
\label{cn}
C_0:=C^2\frac{1-\lambda +(1-\theta+\lambda\theta)C^2}{ 1-\lambda+\lambda\theta C^{-2}},
\feqn
where $C$ is the constant introduced in \eqref{cis}. Then, for all $(\lambda,\theta)\in\Omega^\circ,$
\beqn
\label{apu}
0<E(\ln u_n)-\nu<
C_0\sum_{m=n}^\infty  E\bigg(\prod_{k=0}^m \frac{(1-\lambda)(1-\theta)}{(1-\lambda+\lambda\theta\gamma_k)^2} \bigg).
\feqn
\end{proposition}
The proof of the proposition is included in Section~\ref{pier}. Notice that if the pairs $(\veps_n,\delta_n)$ are independent of each other,
\eqref{apu} yields the error bound
\beq
0<E(\ln u_n)-\nu<\frac{C_0}{1-h}h^{n+1},
\feq
where
\beqn
\label{coh}
h:=\left\{
\begin{array}{ll}
E\Big(\big(\frac{1-\lambda}{1-\lambda+\lambda\theta\gamma_0} \big)^2\Big)&\mbox{\rm if}~\lambda\leq \theta,
\\[3mm]
E\Big(\big(\frac{1-\theta}{1-\theta+\lambda\theta\gamma_0} \big)^2\Big)&\mbox{\rm if}~\lambda\geq \theta.
\end{array}
\right.
\feqn
We remark that an exponential rate of convergence is generic for iteration of random i.\,i.\,d. maps (see, for instance, Theorem~1.1 in \cite{diaf}). The bound in \eqref{apu} implies the exponential
convergence in the general case when the sequence of pairs $(\veps_n,\delta_n)_{n\in\zz}$ is stationary and ergodic, but not necessarily i.\,i.\,d. Indeed, $\frac{(1-\lambda)(1-\theta)}{(1-\lambda+\lambda\theta\gamma_k)^2} \leq h_0,$ where
\beqn
\label{hno}
h_0:=\left\{
\begin{array}{ll}
\big(\frac{1-\lambda}{1-\lambda+\lambda\theta C^{-2}} \big)^2&\mbox{\rm if}~\lambda\leq \theta,
\\[2mm]
\big(\frac{1-\theta}{1-\theta+\lambda\theta C^{-2}} \big)^2&\mbox{\rm if}~\lambda\geq \theta.
\end{array}
\right.
\feqn
Thus, we have
\begin{corollary}
\label{core}
Under Assumption~\ref{asu3}, for all $(\lambda,\theta)\in\Omega^\circ,$
\beq
0<E(\ln u_n)-\nu<\frac{C_0}{1-h_0}h_0^{n+1},
\feq
where constants $C_0$ and $h_0$ are defined, respectively, in \eqref{cn} and \eqref{hno}.
\end{corollary}
The following proposition is an analogue of Proposition~\ref{ester} for the sequence $v_n.$
\begin{proposition}
\label{ester1}
Let Assumption~\ref{asu3} hold and denote
\beq
\witi C_0:=C^2\frac{1-\theta +(1-\lambda+\lambda\theta)C^2}{ 1-\theta+\lambda\theta C^{-2}}.
\feq
Then, for all $(\lambda,\theta)\in\Omega^\circ,$
\beq
0<E(\ln v_n)-\nu<
\witi C_0\sum_{m=n}^\infty  E\bigg(\prod_{k=0}^m \frac{(1-\lambda)(1-\theta)}{(1-\lambda+\lambda\theta\zeta_k)^2} \bigg)
<\frac{\witi C_0}{1-h_0}h_0^{n+1},
\feq
where $h_0$ is defined in \eqref{hno}. If the pairs $(\veps_n,\delta_n)$ are i.\,i.\,d., then the bound can be improved to
\beq
0<E(\ln v_n)-\nu<\frac{\witi C_0}{1-h}h^{n+1},
\feq
where $h$ is defined in \eqref{coh}.
\end{proposition}
We conclude this section with a brief general discussion of Theorem~\ref{verba}. The theorem is closely related to several results that are known for fairly general class of random matrix products in the i.\,i.\,d. case, see for instance \cite{darm} for a brief summary. First, we observe the following. For $n\in\zz,$ let $\calf_n$ and $\calg_n$ denote the $\sigma$-algebras generated by, respectively, $(p_k,q_k)_{k\leq n}$ and $(r_k,s_k)_{k\leq n}.$
\begin{lemma}
\label{station}
Suppose that Assumption~\ref{asu3} is satisfied, and let $u(\mz),$ $v(\mz),$ $u,$ and $v$ be as defined in the statement of Theorem~\ref{verba}.
Then, the following holds true:
\begin{itemize}
\item [(i)] $u(\mz)$ is independent of $\mz$ for $0\leq \mz <(1-\lambda)^{-1}.$
\item [(ii)] $v(\mz)$ is independent of $\mz$ for $0\leq \mz <(1-\theta)^{-1}.$
\item [(iii)] Let $z_0=u$ and define $z_n$ for $n\geq 1$ recursively by the formula
$z_n=p_n+q_n/z_{n-1}.$ Then, $(z_n)_{n\geq 0}$ is a strictly positive, stationary and ergodic sequence. Furthermore, $u$ is the unique distribution of $z_0\in \calf_0$ with this property.
\item [(iv)] Let $t_0=v$ and define $t_n$ for $n\geq 1$ recursively by the formula
$t_n=r_n+s_n/t_{n-1}.$ Then, $(t_n)_{n\geq 0}$ is a strictly positive, stationary and ergodic sequence. Furthermore, $v$ is the unique distribution of $t_0\in \calg_0$ with this property.
\end{itemize}
\end{lemma}
The proof of the lemma is given in Section~\ref{stats}. In view of \eqref{xyrec}, the lemma describes the stationary distributions of the ratios
$ \frac{X_{n+1}}{X_n}$ and $ \frac{Y_{n+1}}{Y_n}$. Thus, Theorem~\ref{verba} implies that
\beqn
\label{lira}
\nu=\lim_{n\to\infty} E\Big(\ln \frac{X_{n+1}}{X_n}\Big)=\lim_{n\to\infty} E\Big(\ln \frac{Y_{n+1}}{Y_n}\Big).
\feqn
We remark that part \emph{(iii)} (respectively, \emph{(iv)}) of the lemma is a consequence of part \emph{(i)} (correspondingly, \emph{(ii)}) combined with a variation of Letac contraction principle identifying stationary distributions of (originally) Markov chains; see, for instance, Proposition~1 in \cite{chat}.
\par
Similar to \eqref{lira} results for a general class of i.\,i.\,d. $2\times 2$ matrices follow, for instance, Theorem~5.1 in \cite{mann} and the results in Appendix~C of \cite{mark} (see also more general and abstract Proposition~3.3 in \cite{boug}, Proposition~4.1 in \cite{fukir}, and Lemma~4.3 in \cite{mogul}). Using \eqref{eqa}, one can express the limiting (and hence, unique stationary) distribution of $\frac{Y_n}{X_n}$ in terms of either $u$ or $v,$ cf. \eqref{ratio} below. The limiting  distribution of $\frac{Y_n}{X_n},$ call it $\psi,$ can be then used to obtain another formal expression for $\nu,$ see Proposition~\ref{nep} below. We remark that $\psi$ is a central object of study of \cite{mann}. Relations between continuous fractions, long-term behavior of matrix products and related invariant measures have been studied by many authors starting with \cite{furst}. In particular, \cite{furst} gives a similar to Theorem~\ref{verba} result for a general i.\,i.\,d. case with $p_nq_n<0,$ see the corollary on pp.~387-388 there.
A brief account of major results for $2\times 2$ i.\,i.,d. matrices can be found, for instance, in the introductory sections of the monograph \cite{boug} and research articles \cite{darm, mann}.
\subsection{Alternative continued fractions representation of \pdfmath{\nu}}
\label{fur}
Let $\rr_+$ denote the set of strictly positive real numbers. For a $2\times 2$ matrix $A=(a_{ij})$ with strictly positive entries define a function (``real-valued M\"{o}bius map", cf.~\cite{chat, brocot}) $T_A:\rr_+\to \rr_+$ by setting
\beq
T_A(x)=\frac{a_{11}x+a_{12}}{a_{21}x+a_{22}}.
\feq
It is easy to verify that $T_{AB}(x)=T_A\big(T_B(x)\big)$ (the group property of M\"{o}bius transformations \cite{beardon}). Using this fact and induction, we obtain (cf. \cite{brocot})
\beqn
\label{tm}
T_{M_n\cdots M_0}(x)=\frac{\theta\delta_n}{1-\theta}+ \frac{(1-\theta)^{-1}\kern6em}{\displaystyle
 \frac{\lambda\veps_n }{1-\lambda}+ \frac{(1-\lambda)^{-1}\kern5em}{\displaystyle
    \frac{\theta\delta_{n-1}}{1-\theta} +\stackunder{}{\ddots\stackunder{}{\displaystyle
      {}+ \frac{(1-\theta)^{-1}}{\displaystyle
        \frac{\lambda\veps_0 }{1-\lambda} + \frac{(1-\lambda)^{-1}}{x}}}
}}}.
\feqn
In particular, assuming that $Y_0=0,$
\beq
T_{M_n\cdots M_0}(+\infty):=\lim_{x\to\infty} T_{M_n\cdots M_0}(x)=\frac{(M_n\cdots M_0)_{11}}{(M_n\cdots M_0)_{21}}=\frac{X_{n+1}}{Y_{n+1}}.
\feq
The ratio $\frac{X_{n+1}}{Y_{n+1}}$ converges in distribution to a non-degenerate random variable under some general conditions (see, for instance, Theorem~4.2 in \cite{mann}). Since by \eqref{emn},
\beqn
\label{ratio}
\frac{X_{n+1}}{X_n}=1-\lambda +\lambda\theta \veps_n\delta_n+ \theta\delta_n\frac{Y_n}{X_n},
\feqn
this convergence would imply that
\beq
1-\lambda +\lambda\theta \gamma_0+ \frac{\theta\delta_0}{T_{M_{-1}\cdots M_{-(n+1)}}(+\infty)}
\feq
converges to the stationary distribution of $\frac{X_{n+1}}{X_n}.$ In our case, the convergence of $\frac{X_{n+1}}{X_n}$ has been established in Theorem~\ref{verba}, and can be used along with \eqref{ratio} to claim the convergence in distribution of $\frac{Y_n}{X_n}.$ Thus, we have proved the following:
\begin{proposition}
\label{nep}
Let Assumption~\ref{asu3} hold. Then, the following holds for all $(\lambda,\theta)\in \Omega^\circ:$
\begin{itemize}
\item [(i)] $T_{M_{-1}\cdots M_{-(n+1)}}(+\infty)$ converges, as $n\to\infty,$ to the stationary distribution of $\frac{Y_n}{X_n}.$
\item [(ii)] Moreover,
\beq
\nu(\lambda,\theta)=\lim_{n\to\infty} E\Big(\ln \Big(1-\lambda +\lambda\theta \gamma_0+ \frac{\theta\delta_0}{T_{M_{-1}\cdots M_{-(n+1)}}(+\infty)}\Big)\Big).
\feq
\end{itemize}
\end{proposition}
Note that under our assumptions, the second part of the proposition follows from part \emph{(i)} by the bounded convergence theorem.
As usual, a counterpart of the proposition for matrices $\witi M_n$ can be obtained by using \eqref{mtilde} and \eqref{tilde} in place of \eqref{eqa}.
\subsection{A special case when \pdfmath{\nu(\lambda,\lambda)} can be calculated}
\label{spec}
In this subsection we calculate $\nu(\lambda,\lambda)$ for a special case when both $\veps_0$ and $\theta_0$ are both gamma-distributed.
The example is a straightforward modification of a part of Theorem~1 in \cite{letac} and Theorem~7.1 in \cite{mann}. If $X$ is a positive absolutely continuous random variable, we write either $X\,{\scriptstyle \overset{rv\to df}{\sim}}\, f$ or $f\,{\scriptstyle \overset{df\to rv}{\sim}}\,X$ to indicate that the distribution of $X$ has density $f.$
\par
Recall that a positive absolutely continuous  random variable $\xi$ is gamma-distributed with parameters $h,a>0$ if its probability density function is
\beq
f_{h,a}(x)=\frac{a^{-h}}{\Gamma(h)}x^{h-1}e^{-\frac{x}{a}}, \qquad x>0,
\feq
where $\Gamma$ is the usual gamma function, that is $\Gamma(h)=\int_0^\infty t^{h-1}e^{-t}\,dt.$ If $h$ is integer and $\xi\,{\scriptstyle \overset{rv\to df}{\sim}}\, f_{h,a},$ then $\xi$ is distributed as a sum of $h$ independent exponential distributions
with mean $a.$ Consider now the generalized inverse Gaussian distribution \cite{bnh} with density
\beq
g_{h,a,b}(x)=\frac{a^{\frac{h}{2}}b^{-\frac{h}{2}}}{2K_{|h|}(\sqrt{ab})}x^{h-1}e^{-\frac{1}{2}\big(ax+\frac{b}{x}\big)},\qquad x>0,
\feq
where $h\in\rr$ and $a,b>0$ are parameters and $K_h$ is the hyperbolic Bessel function of the second kind, that is for $h>0$ we have \cite{watson}:
\beq
K_h(x)=\frac{1}{\sqrt{\pi}}\Gamma\big(h+\frac{1}{2}\big)(2x)^h\int_0^\infty \frac{\cos tdt}{(t^2+x^2)^{h+1/2}}.
\feq
For any constants $h,a,b,c>0,$ we have \cite{letac}:
\beqn
\nonumber
(i)~~&&\qquad
\mbox{\rm if}~\xi\,{\scriptstyle \overset{rv\to df}{\sim}}\,  f_{h,a}~\mbox{\rm and}~ \eta \,{\scriptstyle \overset{rv\to df}{\sim}}\,  g_{h,a,b},\quad \mbox{\rm then}\quad c\xi\,{\scriptstyle \overset{rv\to df}{\sim}}\,  f_{h,ac}~\mbox{\rm and}~c\eta \,{\scriptstyle \overset{rv\to df}{\sim}}\,  g_{h,\frac{a}{c},bc}
\\
\label{pl}
(ii)~&&\qquad
\xi \,{\scriptstyle \overset{rv\to df}{\sim}}\,  g_{h,a,b}\quad \mbox{\rm if and only if}\quad \frac{1}{\xi}\,{\scriptstyle \overset{rv\to df}{\sim}}\,  g_{-h,b,a}
\\
\nonumber
(iii)&&\qquad
g_{h,a,b}=g_{-h,a,b}*f_{h,\frac{2}{a}},
\feqn
where $f*g,$ as usual, denotes the convolution: $f*g(x)=\int_0^x f(t)g(x-t)dt,$ $x>0.$
\par
Suppose now that $\veps_0$ and $\delta_0$ are independent of each other and, furthermore, for some constants $h,a,b,c>0,$ the pair $(\veps_0,\delta_0)$ is independent of a random variable $\xi$ such that
$c(1-\lambda)\xi \,{\scriptstyle \overset{rv\to df}{\sim}}\,  g_{h,a,b},$ $\frac{\lambda\veps_0}{c(1-\lambda)}\,{\scriptstyle \overset{rv\to df}{\sim}}\,  f_{h,\frac{2}{b}},$ and $c\lambda\delta_0\,{\scriptstyle \overset{rv\to df}{\sim}}\,  f_{h,\frac{2}{a}}.$  Then, by virtue of \eqref{pl},
\beq
\frac{\lambda\veps_0}{1-\lambda}+\frac{(1-\lambda)^{-1}}{\xi}\,{\scriptstyle \overset{rv\to df}{\sim}}\,  g_{h,\frac{b}{c},ac},
\feq
and
\beq
c(1-\lambda)\bigg(\frac{\lambda\delta_0}{1-\lambda}+ \frac{(1-\lambda)^{-1}\kern2em}{\displaystyle
 \frac{\lambda\veps_0 }{1-\lambda}+ \frac{(1-\lambda)^{-1}\kern1em}{\displaystyle
    \xi }}\bigg)\,{\scriptstyle \overset{rv\to df}{\sim}}\,  f_{h,\frac{2}{a}}*g_{-h,a,b}=g_{h,a,b}\,{\scriptstyle \overset{df\to rv}{\sim}}\,  c(1-\lambda)\xi.
\feq
Hence, by virtue of \eqref{tm} and Proposition~\ref{nep}, $\xi$ is the stationary distribution of $\frac{Y_n}{X_n}.$ Note that the uniqueness of the latter follows, for instance, from \eqref{ratio} and part \emph{(iii)} of Lemma~\ref{station} (alternatively, and more directly, from \eqref{tm} and the same Letac contraction principle we used to prove Lemma~\ref{station}). Thus, in view of \eqref{ratio}, we have proved the following:
\begin{theorem}
\label{ezra}
Let Assumptions~\ref{indi} hold. Suppose in addition that
\beq
\veps_0\,{\scriptstyle \overset{rv\to df}{\sim}}\,  f_{h,\frac{2}{r}}\qquad \mbox{\rm and} \qquad \delta_0\,{\scriptstyle \overset{rv\to df}{\sim}}\,  g_{h,\frac{2}{s}}
\feq
for some constants $h,r,s>0.$ Then,
\beq
\nu(\lambda,\lambda)=E\big(\ln \big(1-\lambda +\lambda^2 \veps_0\delta_0+ \lambda\delta_0\xi\big)\big),
\feq
where $\xi\,{\scriptstyle \overset{rv\to df}{\sim}}\,  g_{-h,a,b}$ with $a=\frac{r}{\lambda},$  $b= \frac{s(1-\lambda)}{\lambda},$ and  $\xi$ is independent of $(\veps_0,\delta_0).$
\end{theorem}
Note that the distribution of $\xi$ defined in the statement of the theorem depends on the value of the parameter $\lambda.$ In the case when $\veps_0$ and $\delta_0$ are identically distributed, similar random continuous fractions (referred to as random Stieltjes functions) have been considered in \cite{mark} and \cite{pade}. In particular, in the case $r=s,$ an explicit formula for $\nu$ in terms of Bessel functions can be obtained using Theorem~4 of \cite{mark}. In fact, using the fact that the two matrices in factorization \eqref{facts} are i.\,i.\,d. when $r=s$ in the conditions of Theorem~\ref{ezra}, one can show that in that case
\beq
\nu(\lambda,\lambda)=\frac{(1-\lambda)^{-h/2}}{K_{|h|}\big(\frac{s\sqrt{1-\lambda}}{\lambda}\big)}\int_0^\infty x^{h-1}\exp\Big(-\frac{1}{2}\Big(\frac{sx}{\lambda}+\frac{s(1-\lambda)}{\lambda x}\Big)\Big)\ln x \,dx.
\feq
\begin{example}
\label{key}
Consider the setting of Theorem~\ref{ezra} with $r=s=8$ and $h=4.5.$ Then
\beq
\nu(0.69,0.69)\approx 0.062518,
\feq
while, utilizing the results of Section~\ref{boundary}, in particular part (ii) of Proposition~\ref{sp1},
\beq
\max \nu(\lambda,\theta)_{(\lambda,\theta)\in \partial \Omega}=\nu(0.53..,1)=\nu(1,0.53..)\approx 0.061395.
\feq
This provides another example, in addition to the one given in Section~\ref{key1}, when ``Kelly's effect" \eqref{kef} holds true.
\end{example}
\section{Appendix: Proofs}
\label{proofs}
\subsection{Proof of Proposition~\ref{kthm}}
\label{kproof}
Since $\mu(\lambda,\theta)=\mu(\theta,\lambda),$ it is sufficient to prove the monotonicity of $\mu$ in $\lambda.$
It follows from \eqref{mu} that
\beqn
\label{pau}
\frac{\partial \mu}{\partial \lambda}=\frac{1}{2}\Big(\gamma \theta -1 +\frac{(1-\gamma\theta)(\lambda +\theta-\gamma\lambda\theta)+2(\gamma-1)\theta}
{\sqrt{ (\lambda +\theta-\gamma\lambda\theta)^2+4(\gamma-1) \lambda \theta}}\Big).
\feqn
Therefore, $\frac{\partial \mu}{\partial \lambda}=0$ would yield (note that $1-\gamma\theta=0$ along with $\frac{\partial \mu}{\partial \lambda}=0$ would imply $(\gamma-1)\theta=0,$ which is precluded by the conditions of the proposition)
\beq
\Bigl((\lambda +\theta-\gamma\lambda\theta)+2\frac{(\gamma-1)\theta}{1-\gamma\theta}\Bigr)^2
= (\lambda +\theta-\gamma\lambda\theta)^2+4(\gamma-1) \lambda \theta,
\feq
or, equivalently,
\beq
\theta(1-\gamma\theta)(\lambda +\theta-\gamma\lambda\theta)+(\gamma-1)\theta^2
= \lambda \theta(1-\gamma\theta)^2.
\feq
This yields $\gamma\theta(1-\theta)=0,$ which is impossible under the conditions of the proposition. Thus, $\frac{\partial \mu}{\partial \lambda}$
never changes the sign. To complete the proof of the proposition, observe that for any fixed $\theta \in (0,1),$
we have $\lim_{\lambda\to 0}\frac{\partial \mu}{\partial \lambda}=(\gamma-1)$ by virtue of \eqref{pau}. \qed
\subsection{Proof of Proposition~\ref{zth}}
\label{zproof}
Assuming that $\gamma=1,$ and hence $\beta=\alpha^{-1},$ we obtain from \eqref{eqa1} and \eqref{ema} that
\beqn
\label{diff}
\begin{array}{lll}
W_{n+1}-W_n&=&-\lambda(1-\theta)W_n+\theta V_n,
\\
V_{n+1}-V_n&=&\lambda(1-\theta)W_n-\theta V_n,
\end{array}
\feqn
where
\beqn
\label{wen}
W_n:=\alpha U_n.
\feqn
Thus, $W_{n+1}+V_{n+1}=W_n+V_n.$ Hence, at equilibrium, when both the differences at \eqref{diff} are zeroes for all $n\geq 0,$ we must have
\beq
\lambda=\frac{\theta}{1-\theta} \frac{V_0}{W_0}.
\feq
On the other hand, if for some $n\geq 0,$ the inequality $\lambda >\frac{\theta}{1-\theta} \frac{V_n}{W_n}$ holds true, then, taking into account that
$V_{n+1}-V_{n+1}=-\big(W_{n+1}-W_{n+1}\big),$ we obtain:
\beq
\frac{\theta}{1-\theta}\frac{V_{n+1}}{W_{n+1}}&=&
\frac{\theta}{1-\theta}\frac{\lambda(1-\theta)W_n+(1-\theta) V_n}{W_n-\lambda(1-\theta)W_n+\theta V_n}=
\frac{\lambda\theta +\theta \frac{V_n}{W_n}}{1-\lambda(1-\theta)+\theta \frac{V_n}{W_n}}
\\
&<&
\frac{\lambda\theta +\lambda(1-\theta)}{1-\lambda(1-\theta)+\lambda(1-\theta)}=\lambda.
\feq
Similarly, if $\lambda <\frac{\theta}{1-\theta} \frac{V_n}{W_n}$ then $\lambda <\frac{\theta}{1-\theta} \frac{V_{n+1}}{W_{n+1}}.$
By induction, in either case, both the sequences $W_n$ and $V_n$ are monotone, and hence the limits (possibly infinite)
\beq
U=\lim_{n\to\infty} W_n\qquad \mbox{\rm and}\qquad V=\lim_{n\to\infty} V_n
\feq
exist. Taking the limits on both sides of \eqref{diff} then yields $\frac{V}{W}=\lambda \frac{1-\theta}{\theta},$
completing the proof of the proposition in view of \eqref{wen}. \qed
\subsection{Proof of Proposition~\ref{conte}}
\label{pconte}
Let
\beq
\Omega_{-1}=\big\{(\lambda,\theta)\in\Omega: (1-\lambda)(1-\theta)\neq 0\big\}.
\feq
Observe that $\det M_n=(1-\lambda)(1-\theta)$ is strictly positive when $(\lambda,\theta)\in \Omega_{-1}.$
Thus, $M_n$ is an invertible matrix with probability one on $\Omega_{-1}.$
Furthermore, due to \eqref{cis},  $\|M_n\|$  and $\|M_n^{-1}\|$ are bounded on $\Omega_{-1}$ for any sub-multiplicative matrix norm $\|\cdot\|.$
That is, there exists a constant $c>0$ such that
\beq
P(c^{-1}\leq \|M_n\|\leq c)=1.
\feq
This follows, for instance, from the fact that for $(\lambda,\theta)\in \Omega^\circ,$ the spectral norm $\|M_n\|_2$ is
the logarithm of the Perron-Frobenius eigenvalue of $M_n,$ and hence is bounded form below by the minimal and from above by the maximal
entry of the matrix. Moreover, $\|M_n\|_2=1$ when $\lambda\theta=0.$
\par
Hence, $\nu(\lambda,\theta)$ is a continuous function of its parameters on $\Omega_{-1}$ by Theorem~C in \cite{conti}.
It remains to extend this claim to the case when either $\lambda=1$ or $\theta=1.$  Toward this end, for $\eta>0,$ define
\beq
M_{n,\eta}=\begin{pmatrix}
1-\lambda +\lambda(\theta +\eta)\veps_n\delta_n& (\theta+\eta)\delta_n
\\
\lambda (1-\theta) \veps_n&1+\eta-\theta
\end{pmatrix}
.
\feq
Then, by virtue of \eqref{cis},
\beq
\det M_{n,\eta} =(1-\lambda)(1-\eta)+\eta(1-\lambda)+\eta \lambda(\theta+\eta) \veps_n\delta_n>0
\feq
for all $\eta>0$ and $(\lambda,\theta)\in \Omega.$ Let $\nu_\eta(\lambda,\theta)$ be the top Lyapunov exponent of $M_{n,\eta},$ that is
\beq
\nu_\eta=\lim_{n\to\infty}\frac{1}{n}\ln \|M_{n,\eta}\cdots M_{1,\eta}\|=\lim_{n\to\infty}\frac{1}{n}E\big(\ln \|M_{n,\eta}\cdots M_{1,\eta}\|\big).
\feq
By Theorem~C in \cite{conti}, $\nu_\eta(\lambda,\theta)$ is a continuous function on the whole domain $\Omega$ for all
$\eta>0.$ Therefore, using the non-negativity of the matrices under consideration and a standard sub-additivity argument, we obtain that for all $(\lambda_0,\theta_0)\in\Omega,$
\beqn
\nonumber
&& \limsup_{(\lambda,\theta)\to(\lambda_0,\theta_0)} \nu(\lambda,\theta)\leq  \limsup_{(\lambda,\theta)\to(\lambda_0,\theta_0)} \nu_\eta(\theta,\lambda)=
\limsup_{(\lambda,\theta)\to(\lambda_0,\theta_0)} \inf_{n\geq 1}\frac{1}{n}E\big(\ln \|M_{n,\eta}\cdots M_{1,\eta}\|\big)
\\
\label{lims}
&&\qquad
\leq
\inf_{n\geq 1}\frac{1}{n} \limsup_{(\lambda,\theta)\to(\lambda_0,\theta_0)} E\big(\ln \|M_{n,\eta}\cdots M_{1,\eta}\|\big)
 = \nu_\eta(\lambda_0,\theta_0).
\feqn
In the above argument we used Kingman's sub-ergodic theorem to write
\beq
\nu_\eta(\lambda,\theta)=\inf_{n\geq 1}\frac{1}{n}E\big(\ln \|M_{n,\eta}\cdots M_{1,\eta}\|\big)
\feq
and the bounded convergence theorem to establish the convergence of $E\big(\ln \|M_{n,\eta}\cdots M_{1,\eta}\|\big)$
when $(\lambda,\theta)$ approaches $(\lambda_0,\theta_0).$
\par
Next, we observe that when $\theta=1,$ $M_{n,\eta}$ are upper triangular matrices, and therefore \cite{pinkus}
\beqn
\label{maxt}
\nu_\eta(\lambda,1)&=&\max\big\{\ln M_{n,\eta}(1,1),\ln M_{n,\eta}(2,2)\big\}
\\
\nonumber
&=&
\max\big\{E\big(\ln \big(1-\lambda +\lambda(1+\eta)\veps_n\delta_n)\big),\,\ln \eta\}.
\feqn
Since $\eta>0$ is arbitrary, in view of \eqref{lims}, we can conclude that for any $\lambda_0\in [0,1],$
\beqn
\label{upper1}
\limsup_{(\lambda,\theta)\to(\lambda_0,1)} \nu_\eta(\lambda,\theta)\leq E\big(\ln \big(1-\lambda +\lambda\veps_n\delta_n)\big).
\feqn
On the other hand, with $K_n$ defined as
\beq
K_n=\begin{pmatrix}
1-\lambda +\lambda\theta\veps_n\delta_n& \theta\delta_n
\\
0&1-\theta
\end{pmatrix}
,
\feq
we obtain that
\beq
\nu(\theta,\lambda)&\geq&  \lim_{n\to \infty}\frac{1}{n} E\big(\ln \|K_n\cdots K_1\|\big)
\\
&=&
\max\big\{E\big(\ln \big(1-\lambda +\lambda\veps_n\delta_n)\big),\,\ln (1-\theta)\},
\feq
where, similarly to \eqref{maxt}, we used the explicit formula for the Lyapunov exponent of triangular matrices  in terms
of the maximum of the expected values of the logarithm of their diagonal terms. The last formula implies that
\beq
\lim_{(\lambda,\theta)\to(\lambda_0,1)}\nu(\lambda,\theta)\geq
E\big(\ln \big(1-\lambda_0 +\lambda_0\veps_n\delta_n)\big).
\feq
Combining this lower bound with \eqref{upper1}, we conclude that
\beq
\lim_{(\lambda,\theta)\to(\lambda_0,1)}\nu(\lambda,\theta)=
\ln E\big(\big(1-\lambda_0 +\lambda_0\veps_n\delta_n)\big)=\nu(\lambda_0,1),
\feq
where the last identity is another application of the formula for the Lyapunov exponent of triangular matrices.
\par
To conclude the proof, it remains to observe that the case $\lambda=1$ can be handled by using the counterpart for $\theta=1$ and appealing to the symmetry identity stated in Proposition~\ref{duet}. \qed
\subsection{Proof of Proposition~\ref{suc}}
\label{psuc}
{\bf (a)}
For a $2\times 2$ matrix $A,$ let $\|A\|$ denote its spectral norm. If follows from \eqref{limas}, Jensen's inequality, and \eqref{eqa1}
that
\beq
\nu=\lim_{n\to\infty}\frac{1}{n}E (\ln X_n ) \leq \lim_{n\to\infty}\frac{1}{n} \ln E (X_n )=\ln \|M\|,
\feq
where $M=E(M_0)$ as in \eqref{ema}. It is easy to see that if $COV(\veps_n,\delta_n)\leq 0,$ then
\beq
\|M\|\leq \mu,
\feq
where $\mu$ is $\|M\|$ in the case when $COV(\veps_n,\delta_n)=0,$ and is calculated in \eqref{mu}.
Thus, by virtue of \eqref{rez1}, $\nu \leq \ln \mu \leq 0$ when $\gamma\leq 1.$ The completes the prove of part (a) of the proposition.\\
$\mbox{}$
\\
{\bf (b)} The claim follows from Propositions~\ref{sp1} and Proposition~\ref{conte}. \qed
\subsection{Proof of Theorem~\ref{verba}}
\label{pverba}
We will only prove the claim for the sequence $u_n,$ the proof for $v_n$ is similar. Let $z_n=\frac{X_{n+1}}{X_n}$ and recall \eqref{pq}.
It follows from \eqref{eqa} and \eqref{emn}, that using a suitable $Y_0,$ we can set $z_0$ to be any desired positive number given any $X_0>0.$ We will therefore assume, without loss of generality, that $X_0=1$ and $z_0=p_0,$ that is $Y_0=(1-\theta)^{-1}\delta_{-1}^{-1}.$ By ``without loss of generality" we mean that the convergence in \eqref{liml} and \eqref{limasol} remains to occur with probability one. It follows then from \eqref{ar} that $z_{n+1}=p_n+\frac{q_n}{z_n},$ and hence, taking in account that $z_0=p_0,$
\beq
z_{n+1}=p_n + \frac{q_n\kern6em}{\displaystyle
  p_{n-1} + \frac{q_{n-1}\kern5em}{\displaystyle
    p_{n-2} +\stackunder{}{\ddots\stackunder{}{\displaystyle
      {}+ \frac{q_2}{\displaystyle
        p_1 + \frac{q_1}{p_0}}}
}}}
\feq
In view of \eqref{liml} and \eqref{limasol}, taking in account our assumption that $X_0=1,$
\beq
\nu=\lim_{n\to\infty} \frac{1}{n}\sum_{k=1}^n \ln z_k=\lim_{n\to\infty} \frac{1}{n}\sum_{k=1}^n E(\ln z_k)=\lim_{n\to\infty} \frac{1}{n}\sum_{k=1}^n E(\ln u_k),\qquad \as
\feq
where the sequence $u_k$ is defined in \eqref{ubc}, and the last equality holds because $u_k$ and $z_k$ are identically distributed
for any fixed $k\geq 0.$
\par
It is easy to verify that, since $q_n<0$ with probability one, $u_n$ is a monotone decreasing sequence. Therefore,
$u=\lim_{n\to\infty} u_n$ exists with probability one. Furthermore, using the continuity property of probability measures and the monotonicity of the $u_k$ sequence,
\beq
P(u<0)&=&P\big(\cup_{k=1}^\infty \{u_k<0\}\big)=\lim_{n\to\infty} P\big(\cup_{k=1}^n \{u_k<0\}\big)
\\
&=& \lim_{n\to\infty} P(u_k<0) =\lim_{n\to\infty} P(z_k<0)=0.
\feq
Thus, $\ln u$ is well-defined if we use the convention that $\ln 0=-\infty.$ Moreover, by the monotone convergence theorem, $E(\ln u)=\lim_{n\to\infty} E(\ln u_n),$ where the left-hand side in principle can be $-\infty.$ The claim follows now from the Ces\`{a}ro's mean convergence theorem which yields
\beq
\nu=\lim_{n\to\infty} \frac{1}{n}\sum_{k=1}^n E(\ln u_k)=E(\ln u),\qquad \as
\feq
Since we already know that $\nu$ is finite under Assumption~\ref{asu3}, the last identity guarantees that in fact $P(u>0)=1.$ \qed
\subsection{Proof of Proposition~\ref{ester}}
\label{pier}
For $n\in\zz$ and $-(1-\theta)d_{n-1}<x<0,$ define $f_n(x)=\frac{q_n}{p_{n-1}+x}.$ Then,
\beqn
\label{fin}
u_{n+1}(\mz)=p_0+f_0\big(f_{-1}\big(\cdots f_{-(n-1)}(q_{-n}\mz)\cdots\big)\big).
\feqn
Notice that, since $\mz<(1-\lambda)^{-1}$ by the assumptions of the theorem,
\beq
q_{-n}\mz=-(1-\lambda)(1-\theta)d_{-n}\mz>-(1-\theta)d_{-n}.
\feq
Furthermore, if $-(1-\theta)d_{-n}<x\leq 0,$ then $f_{-(n-1)}(x)<0$ and
\beq
f_{-(n-1)}(x)=\frac{q_{-(n-1)}}{p_{-n}+x}>\frac{q_{-(n-1)}}{p_{-n}-(1-\theta)d_{-n}}=
\frac{-(1-\lambda)(1-\theta)d_{-(n-1)}}{1-\lambda+\lambda\theta\gamma_{-n}}>-(1-\theta)d_{-(n-1)}.
\feq
Thus, by the induction on $n,$ the expression in the right-hand side of \eqref{fin} is well-defined. Moreover,
for $-(1-\theta)d_{n-1}<x<0,$ we have
\beq
f_n'(x)=\frac{|q_n|}{(p_{n-1}+x)^2}<\frac{|q_n|}{(p_{n-1}-(1-\theta)d_{n-1})^2}=
\frac{(1-\lambda)(1-\theta)d_n}{(1-\lambda+\lambda\theta\gamma_{n-1})^2}.
\feq
It follows that
\beqn
\nonumber
0&<&u_n(\mz)-u_{n+1}(\mz)<\prod_{k=0}^{n-2} \frac{(1-\lambda)(1-\theta)d_{-k}}{(1-\lambda+\lambda\theta\gamma_{-k-1})^2}\cdot \big| f_{-(n-1)}\big(q_{-n}/p_{-n-1}\big)-f_{-(n-1)}(0)\big|
\\
\label{lisp}
&<&
\frac{q_{-n}}{p_{-n-1}}\prod_{k=0}^{n-1} \frac{(1-\lambda)(1-\theta)d_{-k}}{(1-\lambda+\lambda\theta\gamma_{-k-1})^2}.
\feqn
By virtue of part (ii) of Assumption~\ref{asu3},
\beq
0&<&u_n-u_{n+1}<\big(1-\lambda+\lambda\theta C^2+(1-\theta)C^2\big)\prod_{k=0}^n \frac{(1-\lambda)(1-\theta)d_{-k}}{(1-\lambda+\lambda\theta\gamma_{-k-1})^2}
\\
&<&
\big(1-\lambda+\lambda\theta C^2+(1-\theta)C^2\big)(1-\lambda)^{n+1}(1-\theta)^{n+1}C^2\prod_{k=0}^n \frac{1}{(1-\lambda+\lambda\theta\gamma_{-k-1})^2}.
\feq
Furthermore, it follows from \eqref{fin} and the above induction argument that
\beq
u\geq p_0-(1-\theta)d_0\geq 1-\lambda+\lambda\theta C^{-2},
\feq
and hence, with $C_0$ introduced in \eqref{cn}, we get
\beq
0&<&E(\ln u_n)-E(\ln u)=E\Big(\ln (1+\frac{u_n-u}{u}\Big)\Big)<E\Big(\frac{u_n-u}{u}\Big)
\\
&<&
C_0\sum_{k=n}^\infty  E\Big(\prod_{j=0}^k \frac{(1-\lambda)(1-\theta)}{(1-\lambda+\lambda\theta\gamma_{-j-1})^2} \Big)=
C_0\sum_{k=n}^\infty  E\Big(\prod_{j=0}^k \frac{(1-\lambda)(1-\theta)}{(1-\lambda+\lambda\theta\gamma_j)^2} \Big),
\feq
as desired. \qed
\subsection{Proof of Lemma~\ref{station}}
\label{stats}
We will only prove parts (i) and (iii), parts (ii) and (iv) are proved similarly.
\\
{\bf (i)} Similarly to \eqref{lisp},
\beq
0<|u_n(\mz_1)-u_n(\mz_2)|<\prod_{k=0}^{n-1} \frac{(1-\lambda)(1-\theta)d_{-k}}{(1-\lambda+\lambda\theta\gamma_{-k-1})^2}\cdot | \mz_2-\mz_1|,
\feq
and hence
\beq
\limsup_{n\to\infty} \frac{1}{n}\big|u_n(\mz_1)-u_n(\mz_2)\big|\leq \ln\big((1-\lambda)(1-\theta)\big)-2E\big(\ln(1-\lambda+\lambda\theta\gamma_0)\big)\leq \ln h_0<0,
\feq
where $h_0\in (0,1)$ is a constant defined in \eqref{hno}. Since both sequences $u_n(\mz_1)$ and $u_n(\mz_2)$ converge, as $n\to\infty,$
the limit is common and is equal to $u$ introduced in the statement of Theorem~\ref{verba}.
\\
$\mbox{}$
\\
{\bf (iii)} By induction, it follows from \eqref{ubc} that $z_n\in\calf_n$ and that $(z_n)_{n\geq 0}$ is a stationary sequence.
In fact, the Markov chain argument of Letac (see, for instance, Proposition~1 in \cite{chat}) can be applied verbatim to the ``background Markov chain" $W_n=(\veps_k,\delta_k)_{k\leq n}.$ Since the sequence of pairs $(\veps_n,\delta_n)_{n\in \zz}$ is assumed to be ergodic, $z_n\in\calf_n$ proves the ergodicity of $(Z_n)_{n\geq 0}.$  Finally, again just as in the original Letac contraction principle, the claimed uniqueness property of the stationary distribution follows from the fact that with any initial distribution $z_0$ such that $P(z_0>1-\lambda)=1,$ Theorem~\ref{verba}
insures the weak convergence of $z_n$ to the distribution of $u.$
\qed
\subsection{Numerical algorithm and error bounds for \pdfmath{{\rm BERN}(0.3,0.2;2)} in Section~\ref{key1}}
\label{ebn}
In what follows we adopt some notation of \cite{peres}, that we now proceed to introduce within our specific two-dimensional content.
Suppose that, as in {\rm BERN}$(\veps,\delta;\eta)$ of Assumption~\ref{asu4}, each of the variables $\veps_0$ and $\delta_0$ takes value in a finite set. Consequently, the distribution of matrices $M_n$ is supported on a finite set $\{A_1,\ldots,A_b\}.$ For in instance, $b=4$ in all instances of {\rm BERN}$(\veps,\delta;\eta).$ Denote $p_i=P(M_0=A_i).$ For $x=(x_1,x_2)\in\rr^2,$ let $\|x\|=\max_i |x_i|$ and, for $x\in\rr^2\backslash\{(0,0)\},$ let $\ol x$ denote the direction of $x.$ Let $\pp_+$ be the ``positive quadrant" in the projective space of $(\rr^2,\|\cdot\|),$ that is
\beq
\pp_+=\big\{\ol x: x=(x_1,x_2)\in \rr^2,\,x_i>0~\mbox{\rm for}~i=1,2\big\}.
\feq
In other words, $\pp_+$ is the space of the equivalence classes of non-zero vectors in the positive cone of $(\rr^2,\|\cdot\|),$
where two vectors are considered to be equivalent if they share the same direction. For any practical purpose in this paper, $\ol x$ can be thought as the unit vector $\frac{x}{\|x\|},$ and $\pp_+$ can be identified with the first quadrant arch of the unit circle in $\rr^2.$  With any of our matrices $A_j,$ one can associate an
operator $A_{j \bullet}:\pp_+\to\pp_+$ as follows: $A_{j \bullet} \ol x=\ol {A_j x}.$ Notice that
\beq
(AB)_\bullet \ol x=\ol{(AB) x}=\ol{A(B x)}=\ol{A(\ol {B x})}=A_\bullet (B_\bullet \ol x),
\feq
that is $A_\bullet B_\bullet=(AB)_\bullet.$ Let $\calc_2$ be the space of complex-valued functions on $\pp_+,$ and
introduce a \emph{transfer operator} $T:\calc_2\to\calc_2$ by setting
\beq
(Tf)(\ol x)=E\big(f(M_{0\bullet}\ol x)\big)=\sum_{j=1}^b p_j f(A_{j \bullet} \ol x),\qquad \ol x\in \pp_+.
\feq
It turns out that, with $f_j(\ol x)=\log \|A_j \ol x\|,$
\beqn
\label{cores}
\nu=\lim_{n\to\infty} \sum_{j=1}^b p_j(T^n f_j)(\ol x),\qquad \forall~\ol x\in \pp_+.
\feqn
Indeed,
\beq
&&\sum_{j=1}^b p_j(Tf_j)(\ol x)= \sum_{j=1}^b p_j E\big(f_j(M_{0\bullet}\ol x)\big)=E\big(\log \|M_{1\bullet}M_{0\bullet}\ol x\|\big),
\feq
and, more generally, by induction,
\beq
&&\sum_{j=1}^b p_j(T^n f_j)(\ol x)=E\big(\log \|M_{n\bullet}\cdots M_{0\bullet}\ol x\|\big)=E\big(\log \|(M_n \cdots M_0)_\bullet\ol x\|\big).
\feq
Furthermore (see p.~135 in \cite{peres}), for all $\ol x,\ol y\in\pp_+$ we have
\beq
\Big|\sum_{j=1}^b p_j(T^n f_j)(\ol x)-\sum_{j=1}^b p_j(T^n f_j)(\ol y)\Big|\leq \vartheta^n m h(\ol x,\ol y),
\feq
where
\beqn
\label{er3}
\vartheta=\sum_{j=1}^b p_j\tau(A_j)<1,
\feqn
\beq
m=\max_{j=1,\ldots, b}\,\,\sup_{\ol x,\ol y\in \pp_+} \frac{|f_j(\ol x)-f_j(\ol y)|}{h(\ol x,\ol y)},
\feq
with
\beq
h(\ol x,\ol y)=\log \max_{1\leq i,j\leq 2}\frac{x_iy_j}{x_jy_i}
\feq
and
\beq
\tau(A_j)=\sup\Big\{ \frac{h(A_{j \bullet}\ol x,A_{j \bullet}\ol y)}{h(\ol x,\ol y)}:\ol x,\ol y\in\pp_+,\,\ol x\neq \ol y  \Big\}.
\feq
By Birkhoff's formula (see (8) in \cite{peres}),
\beqn
\label{bt}
\tau(A_j)=\frac{1-\sqrt{\psi_j}}{1+\sqrt{\psi_j}},
\feqn
where
\beqn
\label{bt1}
\psi_j=\min_{i,k,m,n}=\frac{A_j(i,k)A_j(m,n)}{A_j(i,n)A_j(m,k)}.
\feqn
It is not hard to check that $m\leq 1$ (see p.~142 in \cite{peres}), and hence, for any $\ol x\in\pp_+$ and $n\in\nn,$ we have
\beqn
\nonumber
&&\Big|\nu-\sum_{j=1}^b p_j(T^n f_j)(\ol y)\Big|\leq \sum_{k=n}^\infty
\Big|\sum_{j=1}^b p_j(T^{k+1} f_j)(\ol x)-\sum_{j=1}^b p_j(T^k f_j)(\ol x)\Big|
\\
\nonumber
&&
\quad
\leq
\sum_{k=n}^\infty \sum_{j=1}^b p_j
\Big|(T^{k+1} f_j)(\ol x)-(T^k f_j)(\ol x)\Big| \leq
\sum_{k=n}^\infty \sum_{j=1}^b p_j \sum_{\ell=1}^b p_\ell
\Big|(T^k f_j)(A_{\ell \bullet} \ol x)-(T^k f_j)(\ol x)\Big|
\\
\nonumber
&&
\quad
\leq \sum_{k=n}^\infty \sum_{j=1}^b p_j \sum_{\ell=1}^b p_\ell \vartheta^k h(\ol x,A_{\ell \bullet} \ol x)
=
\frac{\vartheta^n}{1-\vartheta}\sum_{\ell=1}^b p_\ell h(\ol x,A_{\ell \bullet} \ol x)
\\
\label{er4}
&&
\quad
=
\frac{\vartheta^n}{1-\vartheta}\sum_{\ell=1}^b p_\ell \log \max_{1\leq i,j\leq 2}\frac{x_i (A_\ell x)_j}{x_j(A_\ell x)_i}.
\feqn
Recall \eqref{emn}. If $\ol x=(1,1),$ then
\beq
&&
\max_{1\leq i,j\leq 2}\frac{x_i (M_n x)_j}{x_j(M_n x)_i}=\max\Big\{\frac{1-\lambda +\theta \delta_n[\lambda \veps_n+1]}{(1-\theta) [\lambda\veps_n+1]},\,\frac{(1-\theta) [\lambda\veps_n+1]} {1-\lambda +\theta \delta_n[\lambda\veps_n+1]}\Big\}
\feq
Moreover, using \eqref{bt1},
\beq
\psi(M_n)=\frac{\lambda \theta \veps_n\delta_n}{1-\lambda +\lambda \theta \veps_n\delta_n}.
\feq
Combining this formulas together with \eqref{er3}, \eqref{bt}, and \eqref{er4}, we obtain that for {\rm BERN}$(0.3,0.2;2),$
\beq
\vartheta=0.6088577,
\feq
and, for 20 iterations,
\beqn
\label{er11}
\Big|\nu-\sum_{j=1}^4 p_j(T^{20} f_j)(\ol y)\Big|\leq  0.00003004159,
\feqn
which, as desired, considerably smaller than the gap  between the values shown in \eqref{ma1} and \eqref{ma3}.

\end{document}